\documentclass[10pt,authoryear]{article}
\usepackage{appendix}
\usepackage{setspace}
\usepackage[font=small,  margin=2cm]{caption}
\usepackage[tbtags]{amsmath}
\usepackage{amsthm,amssymb,amsfonts}
\usepackage{natbib}
\usepackage[normalem]{ulem}
\usepackage{lscape}
\usepackage{subfigure}
\usepackage{makecell}
\usepackage{exscale}
\usepackage{booktabs}
\usepackage{array}
\usepackage{fullpage}
\usepackage{url}
\usepackage{algorithm}
\usepackage{algpseudocode}
\usepackage{bm}
\usepackage{smile}
\usepackage{mathtools}
\usepackage{wrapfig}
\usepackage{lipsum}
\usepackage{threeparttable}
\usepackage{mathrsfs}
\usepackage{relsize}
\usepackage{dsfont}
\usepackage{multirow}
\usepackage{apalike}
\usepackage{chngcntr}
\counterwithout{example}{section}
\usepackage[dvipsnames,svgnames,table]{xcolor}
\usepackage[colorlinks=true,linkcolor=blue,urlcolor=blue,citecolor=blue]{hyperref}
\usepackage{xcolor}
\usepackage[resetlabels]{multibib}
\newcites{sec}{Reference}

\def\tr{{\mathop{\text{\rm Tr}}}}
\def\vec{{\mathop{\text{\rm Vec}}}}
\numberwithin{equation}{section}
\numberwithin{theorem}{section}
\numberwithin{corollary}{section}
\numberwithin{definition}{section}

\begin{document}	
\title{\LARGE Generalized Principal Component Analysis for Large-dimensional Matrix Factor Model}
\author{Yong He\footnotemark[1],~~Yujie Hou\footnotemark[1],~~Haixia Liu\footnotemark[1],~~Yalin Wang\footnotemark[2]}
\renewcommand{\thefootnote}{\fnsymbol{footnote}}
\footnotetext[1]{Corresponding author. Institute for Financial Studies, Shandong University, China; e-mail:{\tt heyong@sdu.edu.cn}}
\footnotetext[2]{Corresponding author. School of Mathematics, Shandong University, China; e-mail:{\tt wangyalin@mail.sdu.edu.cn}}
\maketitle
Matrix factor models  have been growing popular dimension reduction tools for large-dimensional matrix time series. However, the  heteroscedasticity of the idiosyncratic components has barely received any attention. Starting from the pseudo likelihood function, this paper introduces a Generalized Principal Component Analysis (GPCA) method for  matrix factor model which takes the heteroscedasticity into account. Theoretically, we first derive the asymptotic distributions of the GPCA estimators by assuming  the separable covariance matrices are known in advance. {We then propose adaptive thresholding estimators for the  separable covariance matrices and derive their convergence rates, which is of independent interest. We also show that this would not alter the asymptotic distributions of the GPCA estimators under certain regular sparsity conditions in the high-dimensional covariance matrix estimation literature.}
The GPCA estimators are shown to be more efficient than the state-of-the-art methods under certain heteroscedasticity conditions. Thorough numerical studies are conducted to demonstrate the superiority of our method over the existing approaches. Analysis of a financial portfolio dataset illustrates the empirical usefulness of the proposed method.

\vspace{0.5em}
\section{Introduction}
Large datasets with an increasing number of variable records pose unprecedented challenges to traditional statistical analysis  in the era of ``big-data".  Factor models are one of the most important statistical tools for dimension reduction, extracting latent principal factors that contribute to the most variation of the data.  Since the seminal works of \cite{bai2002determining} and \cite{stock2002forecasting}, large-dimensional approximate factor models have received increasing attention. There have been many studies in this area during the last two decades, including but not limited to \cite{bai2003inferential}, \cite{onatski2009testing}, \cite{ahn2013eigenvalue}, \cite{fan2013large}, \cite{bai2012statistical}, \cite{trapani2018randomized}, \cite{ait2017using}, \cite{ait2020high}, \cite{barigozzi2020sequential}, \cite{barigozzi2024tail} and \cite{kong2024staleness}. These works all require the assumption that the fourth moments (or even higher moments) of factors and idiosyncratic errors are finite, which may be constrictive when datasets especially in areas of finance and macroeconomics exhibit heavy-tailed features. In the past few years, there have been flourishing works on relaxing the moment conditions for factor analysis, see for example the endeavors by \cite{yu2019robust}, \cite{chen2021quantile}, \cite{he2022large} and \cite{barigozzi2024tail}.

Matrix-variate data arise when one observes a group of variables structured in a well defined matrix form, and have been frequently observed in various research areas such as finance, signal processing and medical imaging. \cite{wang2019factor} first proposed a bilinear matrix factor model for matrix time series $\{\Xb_t,1\leq t\leq T\}$:
\begin{equation}\label{1.1}
\underbrace{\Xb_t}_{p_1\times p_2}=\underbrace{\Rb}_{p_1\times k_1}\times \underbrace{\Fb_t}_{k_1\times k_2}\times\underbrace{\Cb^\top}_{k_2\times p_2}+\underbrace{\Eb_t}_{p_1\times p_2},
\end{equation}
where $\Rb$ is the $p_1\times k_1$ row factor loading matrix exploiting the variations of $\Xb_{t}$ across the rows, $\Cb$ is the $p_{2}\times k_{2}$ column factor loading matrix reflecting the differences across the columns of $\Xb_{t}$, $\Fb_{t}$ is the common factor matrix for all cells in $\Xb_{t}$, and $\Eb_{t}$ is the idiosyncratic component. Note that the model $(\ref{1.1})$ not only fully utilizes the multidimensional structures of matrices, but also captures the spatial correlations among the elements of the observation matrix in a parsimonious way. 
This model along with its variants have been studied in \cite{chen2020constrained, chen2022modeling,kong2022matrix,yu2022projected,yuan2023two,he2024matrix,zhang2024modeling}. \cite{wang2019factor} proposed estimators of the factor loading matrices and numbers of the row and column factors based on an eigen-analysis of the auto-cross-covariance matrix, extending the theoretical analysis framework of \cite{lam2012factor} to the matrix-variate data setting. \cite{yu2022projected} proposed a projection-based method which improves the estimation efficiency of the row and column loading matrices. \cite{chen2023statistical}  proposed an $\alpha$-PCA method based on eigen-analysis of the mean-covariance weighted average matrix. \cite{he2024matrix} established the equivalence between minimizing the squared loss and the PE method by \cite{yu2022projected} and further proposed a robust method via replacing the squared loss with Huber loss.  Recently, there have also been some works on high-order tensor factor models, see for example, \cite{han2022rank}, \cite{10.1093/jrsssb/qkae001}, \cite{han2020tensor}, \cite{han2024cp}, \cite{lam2021rank}, \cite{chen2022analysis},  \cite{chen2022factor}, \cite{chen2024rank}, \cite{chang2023modelling} and \cite{barigozzi2024tail}.

\begin{figure}[!ht]
  \centering
  \includegraphics[width=12cm, height=4cm]{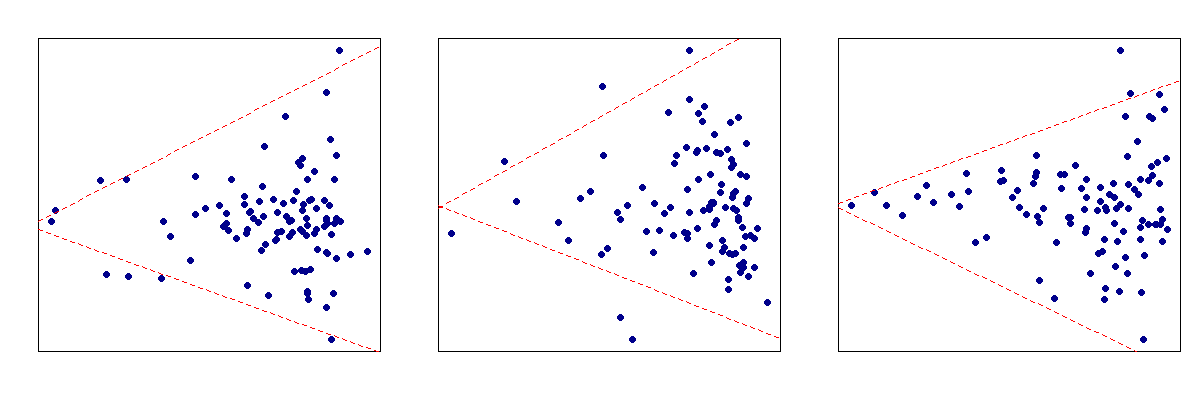}
 \caption{Residual plots  from analyzing a financial portfolio dataset using the PE method proposed by \cite{yu2022projected}. The horizontal axis of each subplot represents the vectorized estimates of  the idiosyncratic component at a certain time point $t$ and the vertical axis represents the corresponding vectorized estimates of the common component.}\label{resi}
 \end{figure}

One drawback of the existing methods for analyzing matrix factor models is that they do not account for the potential heteroscedasticity of the idiosyncratic components.  Figure $\ref{resi}$ illustrates some residual plots from analyzing a financial portfolio dataset using the PE method proposed by \cite{yu2022projected}, from which we can see significant heteroscedasticity of the idiosyncratic components. As is well known, PCA and ordinary least squares (OLS) are equivalent for estimating factor models \citep{fan2011high,he2024matrix}. When heteroscedasticity is present, OLS estimators are not the Best Linear Unbiased Estimators (BLUE). Although OLS estimators remain unbiased, they are not efficient, resulting in less precise parameter estimation. Heteroscedasticity-robust method, weighted least squares (WLS)  is often employed to  address this issue. Few works exist for heteroscedasticity-robust factor analysis to the best of our knowledge. \cite{meijer1996factor} proposed a generalized least squares method based on the first three sample moments in the fixed dimension case. \cite{lewin2003heteroscedastic} proposed a heteroscedastic factor analysis model without specifying the distribution forms for the factors and heteroscedastic errors, and conducted heteroscedastic evaluation and inference of factor scores, but do not provide any theoretical guarantee/analysis. For large-dimensional approximate factor model, \cite{choi2012efficient} derived the conditional maximum likelihood estimators of factors and loadings, along with their asymptotic distributions, which demonstrate that the asymptotic variances of the estimated factors are smaller than those of the corresponding principal component estimators by \cite{bai2003inferential}.
 
For matrix factor model $(\ref{1.1})$, we assume the covariance matrix of vectorized idiosyncratic components has a separable structure, i.e., it can be decomposed as $\Vb\otimes\Ub$, where $\Ub\in\RR^{p_1\times p_1}$ and $\Vb\in\RR^{p_2\times p_2}$ are positive definite matrices with $\Ub$ capturing row-wise dependencies in $\Eb_t$ and $\Vb$ capturing the column-wise dependencies. If $\Ub$ and $\Vb$ are not identity matrices, then the idiosyncratic component $\Eb_t$ exhibits cross-sectional heterogeneity. We propose a \textbf{G}eneralized \textbf{P}rincipal \textbf{C}omponent \textbf{A}nalysis (GPCA) method for   matrix factor model to address the heteroscedasticity of the idiosyncratic components. In the following, we introduce the main idea of the GPCA method briefly. Let $\Zb_t=\Ub^{-1/2}\Xb_t\Vb^{-1/2}$, by reformulating model (\ref{1.1}), we  have
\begin{equation}\label{1.2}
\Zb_t=\Ub^{-1/2}\Rb\Fb_t\Cb^\top\Vb^{-1/2}+\Ub^{-1/2}\Eb_t\Vb^{-1/2}:=\Rb^\ast\Fb_t{\Cb^{\ast}}^\top+\Eb_t^\ast,
\end{equation}
where $\Rb^\ast=\Ub^{-1/2}\Rb$, $\Cb^{\ast}=\Vb^{-1/2}\Cb$, $\Eb_t^\ast=\Ub^{-1/2}\Eb_t\Vb^{-1/2}$. The loading matrices $\Rb^\ast$ and $\Cb^\ast$ in model $(\ref{1.2})$ are not separately identifiable. Without loss of generality, for identifiability, we  further assume that ${\Rb^\ast}^\top\Rb^\ast/p_1=\Ib_{k_1}$ and ${\Cb^\ast}^\top\Cb^\ast/p_2=\Ib_{k_2}$, i.e,  $\Rb^\top\Ub^{-1}\Rb/p_1=\Ib_{k_1}$ and $\Cb^\top\Vb^{-1}\Cb/p_2=\Ib_{k_2}$. We  obtain the maximum pseudo-likelihood estimators of $\Fb_t$, $\Rb$ and $\Cb$, also named as GPCA estimators, by maximizing
\begin{equation}\label{2.1}
2L(\Rb,\Fb_t,\Cb)=-\left\{p_1p_2\ln(2\pi)+p_1\ln|\Vb|+p_2\ln|\Ub|+\tr\left[\Vb^{-1}(\Xb_t-\Rb\Fb_t\Cb^\top)^\top\Ub^{-1}(\Xb_t-\Rb\Fb_t\Cb^\top)\right]\right\}.
\end{equation}
In practice, both $\Ub$ and $\Vb$ are unknown and need to be estimated. We propose an adaptive thresholding estimation method to estimate $\Ub$ and $\Vb$. One point we want to emphasize is that, although the adaptive thresholding estimators and true covariance matrices  $\Ub$ and $\Vb$ differ by constant multiples due to identifiability, this would not bring any trouble for the estimation of the loading spaces and common components by GPCA method. To distinguish, when the  separable covariance matrices are known in advance, we refer to the proposed method as Oracle GPCA, while when $\Ub$ and $\Vb$ are unknown and need to be estimated, it is referred to as data-driven GPCA method.

The contributions of our work lie in the following aspects. Firstly, we for the first time propose the Oracle GPCA method for matrix factor model to address the heteroscedasticity of the idiosyncratic components. Secondly, we derive the asymptotic distributions of the estimated loadings, factors and common components by Oracle GPCA method. Compared with \cite{yu2022projected}, we also provide the asymptotic distributions of the estimated factors and common components by Oracle GPCA method, which is absolutely more challenging than the derivation for vector case in \cite{choi2012efficient} due to the complex matrix structure of  model $(\ref{1.1})$. The asymptotic theoretical analysis indicates that  the GPCA method is more efficient than PE method proposed by \cite{yu2022projected} under certain heteroscedasticity conditions. {In addition, we propose adaptive thresholding estimators for the unknown separable covariance matrices and derive their convergence rates, which is  a much-needed addition to the high-dimensional covariance matrix estimation literature and is of independent interest. At last, we also show that this data-driven GPCA method would not alter the asymptotic distributions under some common sparsity conditions in the high-dimensional covariance matrix estimation literature.}

The rest of this paper is organized as follows. Section $\ref{Section2}$ introduces the Oracle GPCA method  and derives the asymptotic distributions of the estimated loadings, factors and common components by Oracle GPCA. In Section $\ref{Section3}$, we propose adaptive thresholding estimators for the
separable covariance matrices and show that this would not alter the asymptotic distributions of the GPCA
estimators under some sparsity conditions. In Section $\ref{Section4}$, we conduct thorough numerical studies to investigate the finite sample performances of GPCA method over the state-of-the-art methods. Section $\ref{Section5}$ demonstrates the empirical usefulness of the proposed GPCA method by analyzing a financial dataset. We discuss possible future research directions and
conclude the article in Section $\ref{Section6}$. The proofs of the main theorems are collected in the supplementary materials.

To end this section, we introduce some notations that will be used throughout the study. For a real number $a$, let $\textrm{sgn}(a) = 1$ if $a\geq 0$ and $\textrm{sgn} =-1$ if $a<0$. Denote $[n]$ as the set $\{1,\ldots,n\}$ for any positive integer $n$ and ${\rm diag}(a_1,\ldots,a_p)$ as a $p\times p$ diagonal matrix, whose diagonal elements are $a_1,\ldots,a_p$. For a matrix $\Ab$, let $a_{ij}$ or $\Ab_{ij}$ be the $(i,j)$th element of $\Ab$, $\Ab_{i,\cdot}$ be the $i$th row and $\Ab_{\cdot,j}$ be the $j$th column and $\vec(\Ab)$ be the vector obtained by stacking the columns of $\Ab$. Denote $\Ab^\top$ as the transpose of $\Ab$, $\tr(\Ab)$ as the trace of $\Ab$. If $\Ab$ is a symmetric matrix, let $\lambda_j(\Ab)$ be the $j$th largest eigenvalue of $\Ab$; if $\Ab$ is singular, let $\nu_{\min}(\Ab)$ be the minimum singular value of $\Ab$.  $\|\Ab\|$, $\|\Ab\|_F$, $\|\Ab\|_1$ and $\|\Ab\|_{\max}$ represent the spectral norm, Frobenius norm, $L_1$-norm and elementwise norm of $\Ab$, respectively. $\Ib_k$ denotes a $k$-order identity matrix and $\otimes$ denotes the Kronecker product. The notations $\stackrel{p}\rightarrow,\stackrel{d}\rightarrow,\stackrel{a.s.}\rightarrow$ represent convergence in probability, in distribution and almost surely, respectively. The $o_p$ is for convergence to zero in probability and $O_p$ is for stochastic boundedness. For two series of random variables, $X_n$ and $Y_n$, $X_n\lesssim Y_n$ means that $X_n=O_p(Y_n)$, $X_n\gtrsim Y_n$ means that $Y_n=O_p(X_n)$. The notation $X_n\asymp Y_n$ means that $X_n=O_p(Y_n)$ and $Y_n=O_p(X_n)$. The constants $c$ and $C$ are universal and may not be identical in different lines.

\section{Methodology}\label{Section2}
\subsection{Oracle Generalized Principal Component Analysis}
In this section, we first introduce the Oracle GPCA method, and by ``Oracle" we mean that the covariance matrices $\Ub$ and $\Vb$ are known in advance. As stated in the introduction, 
we need to maximize $2L(\Rb,\Fb_t,\Cb)$ in $(\ref{2.1})$ with respect to  $\Rb$, $\Fb_t$ and $\Cb$ to obtain the pseudo likelihood estimators of the factors and loading spaces, which is  equivalent to minimizing the following squared loss under the identifiability conditions
\begin{equation}\label{2.2}
\begin{aligned}
\min_{\{\Rb,\Fb_t,\Cb\}}&L_1(\Rb,\Fb_t,\Cb)=\frac{1}{T}\sum_{t=1}^T\Big\|\Ub^{-1/2}(\Xb_t-\Rb\Fb_t\Cb^\top)\Vb^{-1/2}\Big\|_F^2,\\
&\text{s.t.}\ \frac{1}{p_1}\Rb^\top\Ub^{-1}\Rb=\Ib_{k_1},\quad \frac{1}{p_2}\Cb^\top\Vb^{-1}\Cb=\Ib_{k_2}.
\end{aligned}
\end{equation}
The right hand side of $(\ref{2.2})$ can be simplified as
\begin{equation}\nonumber
\frac{1}{T}\sum_{t=1}^T\Big\|\Ub^{-1/2}(\Xb_t-\Rb\Fb_t\Cb^\top)\Vb^{-1/2}\Big\|_F^2=\frac{1}{T}\sum_{t=1}^T\Big{[}\tr(\Vb^{-1}\Xb_t^\top\Ub^{-1}\Xb_t)-2\tr(\Vb^{-1}\Xb_t^\top\Ub^{-1}\Rb\Fb_t\Cb^\top)+p_1p_2\tr(\Fb_t^\top\Fb_t)\Big{]}.
\end{equation}
The optimization problem is nonconvex jointly in $\{\Rb,\Fb_t,\Cb\}$, but given any two parameters, the loss function is convex over the remaining parameter. For instance, given $\{\Rb,\Cb\}$, $L_1(\Rb,\Fb_t,\Cb)$ is convex over $\Fb_t$. To begin with, we assume that $\{\Rb,\Cb\}$ are given and solve the optimization problem on $\Fb_t$. For each $t$, taking $\partial L_1(\Rb,\Fb_t,\Cb)/\partial\Fb_t=0$, we obtain
\[
\Fb_t=\frac{1}{p_1p_2}\Rb^\top\Ub^{-1}\Xb_t\Vb^{-1}\Cb.
\]
Thus, by substituting $\Fb_t=\Rb^\top\Ub^{-1}\Xb_t\Vb^{-1}\Cb/(p_1p_2)$ into the loss function $L_1(\Fb_t,\Rb,\Cb)$, we further have
\begin{equation}\label{2.3}
\begin{aligned}
&\min_{\{\Rb,\Cb\}}L_1(\Rb,\Cb)=\frac{1}{T}\sum_{t=1}^T\left[\tr(\Vb^{-1}\Xb_t^\top\Ub^{-1}\Xb_t)-\frac{1}{p_1p_2}\tr(\Xb_t^\top\Ub^{-1}\Rb\Rb^\top\Ub^{-1}\Xb_t\Vb^{-1}\Cb\Cb^\top\Vb^{-1})\right],\\
&\quad \quad \text{s.t.}\ \frac{1}{p_1}\Rb^\top\Ub^{-1}\Rb=\Ib_{k_1},\quad \frac{1}{p_2}\Cb^\top\Vb^{-1}\Cb=\Ib_{k_2}.
\end{aligned}
\end{equation}
The Lagrangian function is as follows:
\[
\min_{\{\Rb,\Cb\}}\mathcal{L}_1=L_1(\Rb,\Cb)+\tr\left[\bTheta\left(\frac{1}{p_1}\Rb^\top\Ub^{-1}\Rb-\Ib_{k_1}\right)\right]+\tr\left[\bLambda\left(\frac{1}{p_2}\Cb^\top\Vb^{-1}\Cb-\Ib_{k_2}\right)\right],
\]
where the Lagrangian multipliers $\bTheta$ and $\bLambda$ are symmetric matrices. According to the KKT conditions, let 
\begin{equation}\nonumber
\begin{cases}
\displaystyle\frac{\partial\mathcal{L}_1}{\partial\Rb}=-\frac{1}{T}\sum_{t=1}^T\frac{2}{p_1p_2}\Ub^{-1}\Xb_t\Vb^{-1}\Cb\Cb^\top\Vb^{-1}\Xb_t^\top\Ub^{-1}\Rb+\frac{2}{p_1}\Ub^{-1}\Rb\bTheta=0,\\
\displaystyle\frac{\partial\mathcal{L}_1}{\partial\Cb}=-\frac{1}{T}\sum_{t=1}^T\frac{2}{p_1p_2}\Vb^{-1}\Xb_t^\top\Ub^{-1}\Rb\Rb^\top\Ub^{-1}\Xb_t\Vb^{-1}\Cb+\frac{2}{p_2}\Vb^{-1}\Cb\bLambda=0,
\end{cases}
\end{equation}
respectively, it holds that
\begin{equation}\label{2.4}
\left\{
\begin{aligned}
\left(\frac{1}{Tp_2}\sum_{t=1}^T\Ub^{-1/2}\Xb_t\Vb^{-1}\Cb\Cb^\top\Vb^{-1}\Xb_t^\top\Ub^{-1/2}\right)\Ub^{-1/2}\Rb&=\Ub^{-1/2}\Rb\bTheta,\\
\left(\frac{1}{Tp_1}\sum_{t=1}^T\Vb^{-1/2}\Xb_t^\top\Ub^{-1}\Rb\Rb^\top\Ub^{-1}\Xb_t\Vb^{-1/2}\right)\Vb^{-1/2}\Cb&=\Vb^{-1/2}\Cb\bLambda,
\end{aligned}
\right.
\quad\textrm{or}\quad
\left\{
\begin{aligned}
\Mb_{\Cb}\Ub^{-1/2}\Rb&=\Ub^{-1/2}\Rb\bTheta,\\
\Mb_{\Rb}\Vb^{-1/2}\Cb&=\Vb^{-1/2}\Cb\bLambda,
\end{aligned}
\right.
\end{equation}
where
\[
\Mb_{\Cb}=\frac{1}{Tp_2}\sum_{t=1}^T\Ub^{-1/2}\Xb_t\Vb^{-1}\Cb\Cb^\top\Vb^{-1}\Xb_t^\top\Ub^{-1/2},\quad\Mb_{\Rb}=\frac{1}{Tp_1}\sum_{t=1}^T\Vb^{-1/2}\Xb_t^\top\Ub^{-1}\Rb\Rb^\top\Ub^{-1}\Xb_t\Vb^{-1/2}.
\]
We denote the first $k_1$ eigenvectors of $\Mb_{\Cb}$ as $\{\br(1),\ldots,\br(k_1)\}$ and the corresponding eigenvalues as $\{\theta_1,\ldots,\theta_{k_1}\}$. Similarly, denote the first $k_2$ eigenvectors of $\Mb_{\Rb}$ as $\{\bc(1),\ldots,\bc(k_2)\}$ and the corresponding eigenvalues as $\{\lambda_1,\ldots,\lambda_{k_2}\}$. From $(\ref{2.4})$, we can see that $\Rb=\sqrt{p_1}\Ub^{1/2}(\br(1),\ldots,\br(k_1))$, $\Cb=\sqrt{p_2}\Vb^{1/2}(\bc(1),$ $\ldots,\bc(k_2))$, $\bTheta=\textrm{diag}(\theta_1,\ldots,\theta_{k_1})$, $\bLambda=\textrm{diag}(\lambda_1,\ldots,\lambda_{k_2})$ satisfy the KKT conditions. However, $\Mb_{\Cb}$ depends on the unknown column factor loading $\Cb$ while $\Mb_{\Rb}$ depends on the unknown row factor loading $\Rb$, which prompts us to consider an iterative procedure to obtain the estimators. We summarize the Oracle GPCA procedure in Algorithm $\ref{algori2.1}$. The initialization estimators $\widehat{\Rb}_0$ and $\widehat{\Cb}_0$ could be obtained by setting $\widehat{\Rb}_0=\Ub^{1/2}\widehat{\Rb}_0^\ast$ and $\widehat{\Cb}_0=\Vb^{1/2}\widehat{\Cb}_0^\ast$, where $\widehat{\Rb}_0^\ast$ and $\widehat{\Cb}_0^\ast$ can be any estimators in the literature, depending on the model assumptions, see for example \cite{wang2019factor}, \cite{chen2023statistical}, \cite{he2024matrix}. 

\begin{algorithm}
\setstretch{1.5}
  \caption{Oracle Generalized Principal Component Analysis }
  \label{algori2.1}
  \begin{algorithmic}[1]
    \Require Data matrices $\{\Xb_t\}_{1\leq t\leq T}$, covariance matrices $\Ub$ and $\Vb$, the pair of factor numbers $(k_1,k_2)$, the initial estimators $\widehat{\Rb}_0^\ast$ and $\widehat{\Cb}_0^\ast$
    \Ensure Factor loading matrices $\hat{\Rb}$ and $\hat{\Cb}$
    \State define $\Zb_t=\Ub^{-1/2}\Xb_t\Vb^{-1/2}$;
    \State define $\widehat{\Mb}_{\Cb}=(Tp_2)^{-1}\sum_{t=1}^T\Zb_t\widehat{\Cb}_0^\ast\widehat{\Cb}_0^{\ast\top}\Zb_t^\top$ and $\widehat{\Mb}_{\Rb}=(Tp_1)^{-1}\sum_{t=1}^T\Zb_t^\top\widehat{\Rb}_0^\ast\widehat{\Rb}_0^{\ast\top}\Zb_t$,  obtain the leading $k_1$ eigenvectors of $\widehat{\Mb}_{\Cb}$  and the leading $k_2$ eigenvectors of $\widehat{\Mb}_{\Rb}$, denoted as $\hat{\Qb}_\Rb$ and $\hat{\Qb}_\Cb$, respectively;
    \State the row and column loading matrices are finally given by $\hat{\Rb}=\sqrt{p_1}\Ub^{1/2}\hat{\Qb}_\Rb$ and $\hat{\Cb}=\sqrt{p_2}\Vb^{1/2}\hat{\Qb}_\Cb$.
  \end{algorithmic}
\end{algorithm}
\subsection{Theoretical Results}
In this section, we first derive the asymptotic distributions of the one-step estimators of loadings.
To facilitate the comparison of  theoretical properties of the Oracle GPCA estimators with those of the PE estimators, we select the estimators that satisfy the sufficient condition in \cite{yu2022projected}, such as the $\alpha$-PCA estimators, as the initial estimators $\widehat{\Rb}_0^\ast$ and $\widehat{\Cb}_0^\ast$ in Algorithm \ref{algori2.1}. To derive the asymptotic distributions, we need the following assumptions.

\begin{asmp}[Strong Mixing]\label{asmpa-mixing}
Let $\mathcal{F}_{-\infty}^0$ and $\mathcal{F}_T^{\infty}$ denote the $\sigma$-algebras generated by $\{(\vec(\Fb_t),\vec(\Eb_t^\ast)):t\leq 0\}$ and $\{(\vec(\Fb_t),\vec(\Eb_t^\ast)):t\geq T\}$, respectively. There exist positive constants $r_1$ and $C$ such that for all $T\in\ZZ^{+}$,
$$
\alpha(T)\leq\exp(-CT^{r_1}),
$$
where $\displaystyle\alpha(T)=\sup_{A\in\mathcal{F}_{-\infty}^0,B\in\mathcal{F}_{T}^\infty}\vert P(A\cap B)-P(A)P(B)\vert$.
\end{asmp}
\begin{asmp}[Factor Matrix]\label{asmpFactor}
The factor matrix satisfies $\EE(\Fb_t)=\textbf{0},\EE\|\Fb_t\|^4\leq c<\infty$ for some constant $c>0$, and
\begin{equation}\label{2.5}
\frac{1}{T}\sum_{t=1}^T\Fb_t\Fb_t^\top\stackrel{p}\rightarrow\bSigma_1\quad\textrm{and}\quad\frac{1}{T}\sum_{t=1}^T\Fb_t^\top\Fb_t\stackrel{p}\rightarrow\bSigma_2,
\end{equation}
where $\bSigma_i$ is a $k_i\times k_i$ positive definite matrix with distinct eigenvalues and spectral decomposition $\bSigma_i=\bGamma_i\bLambda_i\bGamma_i^\top,i=1,2$. The factor numbers $k_1$ and $k_2$ are fixed as $\min\{T,p_1,p_2\}\rightarrow\infty$.
\end{asmp}
\begin{asmp}[Loading Matrix]\label{asmpLoading}
Positive constants $\bar{r}$ and $\bar{c}$ exist such that $\|\Rb^\ast\|_{\max}\leq\bar{r},\|\Cb^\ast\|_{\max}\leq\bar{c}$. As $\min\{p_1,p_2\}\rightarrow\infty,\big\|p_1^{-1}\Rb^{\ast\top}\Rb^\ast-\Ib_{k_1}\big\|\rightarrow0$ and $\big\|p_2^{-1}\Cb^{\ast\top}\Cb^\ast-\Ib_{k_2}\big\|\rightarrow0$.
\end{asmp}

{The strong mixing condition in Assumption $\ref{asmpa-mixing}$ allows weak temporal correlations for both the vectorized factor process $\{\vec(\Fb_t)\}$ and the vectorized noise process $\{\vec(\Eb_t^\ast)\}$ and one can easily verify that the $\alpha$-mixing condition in \cite{yu2022projected} is satisfied under the strong mixing condition.} In Assumption $\ref{asmpFactor}$, the factor matrix is centered with bounded fourth moment. According to Corollary 16.2.4 in \cite{athreya2006measure}, the condition ($\ref{2.5}$) is easily satisfied under the strong mixing assumption. The eigenvalues of $\bSigma_i$’s are assumed to be distinct such that the corresponding eigenvectors are identifiable. Assumption $\ref{asmpLoading}$ is strong factor condition, which means that the row and column factors are pervasive along both dimensions. For identifiability, we assume $\big\|p_1^{-1}\Rb^{\ast\top}\Rb^\ast-\Ib_{k_1}\big\|\rightarrow0$ and $\big\|p_2^{-1}\Cb^{\ast\top}\Cb^{\ast}-\Ib_{k_2}\big\|\rightarrow0$ as $\min\{p_1,p_2\}\rightarrow\infty$. Assumptions $\ref{asmpFactor}$ and $\ref{asmpLoading}$ are standard and common in  related literature and similar assumptions are adopted by \cite{yu2022projected}, \cite{chen2023statistical}, \cite{he2023one} and \cite{he2024matrix}.
\begin{asmp}[Weak Correlation of Noise $\Eb_t^\ast$ across Column, Row, and Time]\label{asmpNoise}
A positive constant $c<\infty$ exists such that

1. $\EE e_{t,ij}^\ast=0,\ \EE(e_{t,ij}^{\ast})^8\leq c$\ ;

2. for any $t\in[T],i\in[p_1],j\in[p_2]$,
\begin{equation}\nonumber
(1).\sum_{s=1}^T\sum_{l=1}^{p_1}\sum_{h=1}^{p_2}\vert\EE e_{t,ij}^\ast e_{s,lh}^\ast\vert\leq c,\quad(2).\sum_{l=1}^{p_1}\sum_{h=1}^{p_2}\vert\EE e_{t,lj}^\ast e_{t,ih}^\ast\vert\leq c \ ;
\end{equation}

3. for any $t\in[T],i,l_1\in[p_1],j,h_1\in[p_2]$,
\begin{equation}\nonumber
\begin{aligned}
&(1).\sum_{s=1}^T\sum_{l_2=1}^{p_1}\sum_{h=1}^{p_2}\Big{\vert}\textrm{Cov}(e_{t,ij}^\ast e_{t,l_1j}^\ast,e_{s,ih}^\ast e_{s,l_2h}^\ast)\Big{\vert}\leq c,\quad\sum_{s=1}^T\sum_{l=1}^{p_1}\sum_{h_2=1}^{p_2}\Big{\vert}\textrm{Cov}(e_{t,ij}^\ast e_{t,ih_1}^\ast,e_{s,lj}^\ast e_{s,lh_2}^\ast)\Big{\vert}\leq c,\\
&(2).\sum_{s=1}^T\sum_{l_2=1}^{p_1}\sum_{h_2=1}^{p_2}\left(\Big{\vert}\textrm{Cov}(e_{t,ij}^\ast e_{t,l_1h_1}^\ast,e_{s,ij}^\ast e_{s,l_2h_2}^\ast)\Big{\vert}+\Big{\vert}\textrm{Cov}(e_{t,l_1j}^\ast e_{t,ih_1}^\ast,e_{s,l_2j}^\ast e_{s,ih_2}^\ast)\Big{\vert}\right)\leq c \ ;
\end{aligned}
\end{equation}
where $e_{t,ij}^\ast$ is the $(i,j)$th element of $\Eb_t^\ast$.
\end{asmp}
Assumption $\ref{asmpNoise}$ is essentially an extension of Assumption C in \cite{bai2003inferential} to the matrix regime. Similar conditions can be found in \cite{yu2022projected}, \cite{chen2023statistical} and \cite{barigozzi2024tail}. Assumption $\ref{asmpNoise}$ 2.(1) allows weak correlation of the noises across time, row and column. Assumption $\ref{asmpNoise}$ 2.(2) further controls the column-wise and row-wise correlation of the noises. The correlation of noises up to the second order is controlled by Assumption $\ref{asmpNoise}$ 3.(1). Assumption $\ref{asmpNoise}$ 3.(2) is similar to $\ref{asmpNoise}$ 3.(1), but on different combinations of noise pairs.
\begin{asmp}[Weak Dependence between Factors $\Fb_t$ and Noises $\Eb_t^\ast$]\label{asmpFactorNoise}
There exists a constant $c>0$, such that,

1. for any deterministic vectors $\bv$ and $\bw$ satisfying $\|\bv\|=1$ and $\|\bw\|=1$,
\begin{equation}\nonumber
\EE\left\|\frac{1}{\sqrt{T}}\sum_{t=1}^T(\Fb_t\bv^\top\Eb_t^\ast\bw)\right\|^2\leq c\ ;
\end{equation}

2. for any $i,l_1\in[p_1]$ and $j,h_1\in[p_2]$,
\begin{equation}\nonumber
\begin{aligned}
&(1).\left\|\sum_{h=1}^{p_2}\EE(\bar{\bzeta}_{ij}\otimes\bar{\bzeta}_{ih})\right\|_{\max}\leq c,\quad\left\|\sum_{l=1}^{p_1}\EE(\bar{\bzeta}_{ij}\otimes\bar{\bzeta}_{lj})\right\|_{\max}\leq c,\\
&(2).\left\|\sum_{l=1}^{p_1}\sum_{h_2=1}^{p_2}\textrm{Cov}(\bar{\bzeta}_{ij}\otimes\bar{\bzeta}_{ih_1},\bar{\bzeta}_{lj}\otimes\bar{\bzeta}_{lh_2})\right\|_{\max}\leq c,\quad\left\|\sum_{l_2=1}^{p_1}\sum_{h=1}^{p_2}\textrm{Cov}(\bar{\bzeta}_{ij}\otimes\bar{\bzeta}_{l_1j},\bar{\bzeta}_{ih}\otimes\bar{\bzeta}_{l_2h})\right\|_{\max}\leq c,
\end{aligned}
\end{equation}
where $\bar{\bzeta}_{ij}=\vec(\sum_{t=1}^T\Fb_te_{t,ij}^\ast/\sqrt{T})$.
\end{asmp}
Assumption $\ref{asmpFactorNoise}$ 1 is summarized from the Assumptions F and G.2 in \cite{chen2023statistical} and is also adopted by \cite{yu2022projected} and \cite{he2024online}. Assumption $\ref{asmpFactorNoise}$ 1 simply implies that $\EE(\Fb_t\bv^\top\Eb_t^\ast\bw)\approx\textbf{0}$ and the temporal dependency of the series $\{\Fb_t\bv^\top\Eb_t^\ast\bw\}$ is also weak. Assumption $\ref{asmpFactorNoise}$ 2 controls higher-order correlations between the factor and noise series, where $\bar{\bzeta}_{ij}$ can be viewed as random vectors with fixed dimension and bounded marginal variances (under Assumption $\ref{asmpFactorNoise}$ 1).
\begin{asmp}[Covariance Matrices $\Ub$ and $\Vb$]\label{asmpUV}
There exist constants $c_1,c_2>0$ such that $\lambda_{\min}(\Ub),\lambda_{\min}(\Vb)>c_1$, $\lambda_{\max}(\Ub),\lambda_{\max}(\Vb)<c_2$.
\end{asmp}
\begin{asmp}\label{asmpVariance}
For $i\in[p_1]$,
\begin{equation}\nonumber
\frac{1}{\sqrt{Tp_2}}\sum_{t=1}^T\Fb_t\Cb^\top\Vb^{-1}\be_{t,i\cdot}\stackrel{d}\rightarrow\mathcal{N}(\textbf{0},\Ab_{1i}),\quad\textrm{where}\quad\Ab_{1i}=\lim_{T,p_1,p_2\rightarrow\infty}\frac{1}{Tp_2}\sum_{t=1}^T\EE\left(\Fb_t\Cb^\top\Vb^{-1}\textrm{Cov}(\be_{t,i\cdot})\Vb^{-1}\Cb\Fb_t^\top\right),
\end{equation}
and for $j\in[p_2]$,
\begin{equation}\nonumber
\frac{1}{\sqrt{Tp_1}}\sum_{t=1}^T\Fb_t^\top\Rb^\top\Ub^{-1}\be_{t,\cdot j}\stackrel{d}\rightarrow\mathcal{N}(\textbf{0},\Ab_{2j}),\quad\textrm{where}\quad\Ab_{2j}=\lim_{T,p_1,p_2\rightarrow\infty}\frac{1}{Tp_1}\sum_{t=1}^T\EE\left(\Fb_t^\top\Rb^\top\Ub^{-1}\textrm{Cov}(\be_{t,\cdot j})\Ub^{-1}\Rb\Fb_t\right),
\end{equation}
where $\Ab_{1i}$ and $\Ab_{2j}$ are positive definite matrices whose eigenvalues are bounded away from $0$ and infinity, and $\be_{t,i\cdot}$ and $\be_{t,\cdot j}$ are the $i$th row and $j$th column of $\Eb_t$, respectively.
\end{asmp}
Assumption $\ref{asmpUV}$ requires that $\Ub$ and $\Vb$ are well-conditioned. Assumption $\ref{asmpVariance}$ can be verified through the martingale central limit theorem. It is easy to fulfill under the imposed strong mixing condition and weak correlation assumptions. One can refer to Chapter 16 of \cite{athreya2006measure} for more details. Similar assumptions are imposed in \cite{bai2003inferential}, \cite{yu2022projected}, \cite{chen2023statistical} and \cite{he2023one}. The following theorem  establishes the asymptotic distributions of the one-step estimators of the loadings.

\begin{theorem}\label{theo3.1}
Under Assumptions $\ref{asmpa-mixing}$ to $\ref{asmpVariance}$ with $(k_1, k_2)$ fixed and given, as $\min\{T,p_1,p_2\}\rightarrow\infty$, there exist matrices $\widehat{\Hb}_{\Rb}$ and $\widehat{\Hb}_{\Cb}$, satisfying $\widehat{\Hb}_{\Rb}^\top\widehat{\Hb}_{\Rb}\stackrel{p}\rightarrow\Ib_{k_1}$ and $\widehat{\Hb}_{\Cb}^\top\widehat{\Hb}_{\Cb}\stackrel{p}\rightarrow\Ib_{k_2}$, such that

1. for $i\in[p_1]$,
\begin{equation}\nonumber
\begin{cases}
\sqrt{Tp_2}(\widehat{\Rb}_{i\cdot}-\widehat{\Hb}_{\Rb}^\top\Rb_{i\cdot})\stackrel{d}\rightarrow\mathcal{N}(\textbf{0},\bLambda_1^{-1}\bGamma_1^\top\Ab_{1i}\bGamma_1\bLambda_1^{-1}),\quad\textrm{if}\quad Tp_2=o(\min\{T^2p_1^2,p_1^2p_2^2\}),\\
\displaystyle\widehat{\Rb}_{i\cdot}-\widehat{\Hb}_{\Rb}^\top\Rb_{i\cdot}=O_p\left(\frac{1}{Tp_1}+\frac{1}{p_1p_2}\right),\quad\quad\quad\quad\quad\quad\quad\textrm{if}\quad Tp_2\gtrsim\min\{T^2p_1^2,p_1^2p_2^2\};
\end{cases}
\end{equation}

2. for $j\in[p_2]$,
\begin{equation}\nonumber
\begin{cases}
\sqrt{Tp_1}(\widehat{\Cb}_{j\cdot}-\widehat{\Hb}_{\Cb}^\top\Cb_{j\cdot})\stackrel{d}\rightarrow\mathcal{N}(\textbf{0},\bLambda_2^{-1}\bGamma_2^\top\Ab_{2j}\bGamma_2\bLambda_2^{-1}),\quad\textrm{if}\quad Tp_1=o(\min\{T^2p_2^2,p_1^2p_2^2\}),\\
\displaystyle\widehat{\Cb}_{j\cdot}-\widehat{\Hb}_{\Cb}^\top\Cb_{j\cdot}=O_p\left(\frac{1}{Tp_2}+\frac{1}{p_1p_2}\right),\quad\quad\quad\quad\quad\quad\quad\textrm{if}\quad Tp_1\gtrsim\min\{T^2p_2^2,p_1^2p_2^2\},
\end{cases}
\end{equation}
where $\Rb_{i\cdot}$, $\widehat{\Rb}_{i\cdot}$, $\Cb_{j\cdot}$ and $\widehat{\Cb}_{j\cdot}$ are the $i$th/$j$th row of $\Rb$, $\widehat{\Rb}$, $\Cb$ and $\widehat{\Cb}$,  $\widehat{\Rb}$ and $\widehat{\Cb}$ are the one-step estimators of $\Rb$ and $\Cb$, respectively.
\end{theorem}
Theorem $\ref{theo3.1}$ shows that $\widehat{\Rb}_{i\cdot}$ and $\widehat{\Cb}_{j\cdot}$ consistently estimate the unobservable loadings $\widehat{\Hb}_{\Rb}^\top\Rb_{i\cdot}$ and $\widehat{\Hb}_{\Cb}^\top\Cb_{j\cdot}$ for $i\in[p_1]$ and $j\in[p_2]$, respectively. Note that the rotational matrices $\widehat{\Hb}_{\Rb}$ and $\widehat{\Hb}_{\Cb}$ are different from those of the PE estimators  in \cite{yu2022projected}. By comparison with Theorem 3.2 of \cite{yu2022projected}, we conclude that the convergence rates of the  Oracle GPCA estimators and PE estimators are the same while the asymptotic distributions are different. In fact, it is hard to compare the efficiency of the Oracle GPCA estimators and the PE estimators due to  the potential complex structures of $\Ub$ and $\Vb$. It is obvious that they have the same asymptotic covariance matrix when both $\Ub$ and $\Vb$ are identity matrices. To illustrate the efficiency gain of the Oracle GPCA estimators over the PE estimators, we focus on the  row loading matrix $\Rb$ and  consider a  particular form of $\Ub$ and $\Vb$ in the following example.
\begin{example}
According to Theorem 3.2 in  \cite{yu2022projected}, 
$$\sqrt{Tp_2}(\hat{\Rb}_{i\cdot}-\hat{\Hb}_1^\top\Rb_{i\cdot})\stackrel{d}\rightarrow\mathcal{N}\left(\bm{0},\bLambda_1^{-1}\bGamma_1^\top\Vb_{1i}\bGamma_1\bLambda_1^{-1}\right),$$
where $\bLambda_1$ and $\bGamma_1$ are the same as defined in Assumption $\ref{asmpFactor}$, and $\displaystyle\Vb_{1i}=\lim_{T,p_1,p_2\rightarrow\infty}\frac{1}{Tp_2}\sum_{t=1}^T\EE\left(\Fb_t\Cb^\top\textrm{Cov}(\be_{t,i\cdot})\Cb\Fb_t^\top\right)$.
By the fact that the covariance matrix of $\vec(\Eb_t)$ is $\Vb\otimes\Ub$, we can obtain that $\textrm{Cov}(\be_{t,i\cdot})=u_{ii}\Vb$, where $u_{ii}$ is the $(i,i)$-th element of $\Ub$. Thus, we have
\begin{equation}\nonumber
\Vb_{1i}=\lim_{T,p_1,p_2\rightarrow\infty}\frac{1}{Tp_2}\sum_{t=1}^T\EE\left(u_{ii}\Fb_t\Cb^\top\Vb\Cb\Fb_t^\top\right)\quad\textrm{and}\quad\Ab_{1i}=\lim_{T,p_1,p_2\rightarrow\infty}\frac{1}{Tp_2}\sum_{t=1}^T\EE\left(u_{ii}\Fb_t\Cb^\top\Vb^{-1}\Cb\Fb_t^\top\right).
\end{equation}
Assuming that $\Ub$ and $\Vb$ are both diagonal positive definite matrices. As long as all diagonal elements of $\Vb$ are greater than 1, we have that $\bLambda_1^{-1}\bGamma_1^\top(\Vb_{1i}-\Ab_{1i})\bGamma_1\bLambda_1^{-1}$ is a positive definite matrix, and consequently, the asymptotic variances of the Oracle GPCA estimators are smaller than those of the PE estimators, i.e., the Oracle GPCA estimators are more efficient than the PE estimators.
\end{example}

The factor matrix $\Fb_t$ can be estimated easily by $\widehat{\Fb}_t=\widehat{\Rb}^\top\Ub^{-1}\Xb_t\Vb^{-1}\widehat{\Cb}/(p_1p_2)$. Then the common component is given by $\widehat{\textbf{S}}_t=\widehat{\Rb}\widehat{\Fb}_t\widehat{\Cb}^\top$. To further study the element-wise asymptotic distributions of the one-step estimators of the factors and common components, we need the following assumption.
\begin{asmp}\label{asmpVariance1}$\ $

1. $\displaystyle\frac{1}{p_1p_2}(\Cb\otimes\Rb)^\top(\Vb\otimes\Ub)^{-1}(\Cb\otimes\Rb)\rightarrow\bSigma^{\ast}$, where $\bSigma^{\ast}$ is a positive-definite matrix.

2. For each $i\in[p_1]$, $j\in[p_2]$ and $t\in[T]$, as $T$, $p_1$ and $p_2$ go to infinity simultaneously,
\begin{equation}\nonumber
\begin{pmatrix}
\displaystyle\frac{1}{\sqrt{Tp_2}}\sum_{s=1}^T\Fb_s\Cb^\top\Vb^{-1}\be_{s,i\cdot}\\
\displaystyle\frac{1}{\sqrt{Tp_1}}\sum_{s=1}^T\Fb_s^\top\Rb^\top\Ub^{-1}\be_{s,\cdot j}\\
\displaystyle\frac{1}{p_1p_2}(\Cb\otimes\Rb)^\top(\Vb\otimes\Ub)^{-1}\Vec(\Eb_t)
\end{pmatrix}\stackrel{d}\rightarrow\mathcal{N}(\textbf{0},\bOmega),\quad\textrm{where}\quad\bOmega=\begin{pmatrix}
\Ab_{1i}&\textbf{0}&\textbf{0}\\
\textbf{0}&\Ab_{2j}&\textbf{0}\\
\textbf{0}&\textbf{0}&\bSigma^\ast
\end{pmatrix}.
\end{equation}
\end{asmp}
Assumption $\ref{asmpVariance1}$ 1 is of standard nature, extending  Assumption 2 (ii) in \cite{choi2012efficient} to the matrix regime. Assumption $\ref{asmpVariance1}$ 2 states that joint asymptotic normality holds for the three vectors $\displaystyle\frac{1}{\sqrt{Tp_2}}\sum_{s=1}^T\Fb_s\Cb^\top\Vb^{-1}\be_{s,i\cdot}$, $\displaystyle\frac{1}{\sqrt{Tp_1}}\sum_{s=1}^T\Fb_s^\top\Rb^\top\Ub^{-1}\be_{s,\cdot j}$ and $\displaystyle\frac{1}{p_1p_2}(\Cb\otimes\Rb)^\top(\Vb\otimes\Ub)^{-1}\Vec(\Eb_t)$, and further requires they are asymptotically uncorrelated. The following Theorems $\ref{theo3.2}$ and $\ref{theo3.3}$ establish the element-wise asymptotic normality of the one-step estimators of the factors and common components, respectively.
\begin{theorem}\label{theo3.2}
Under Assumptions $\ref{asmpa-mixing}$ to $\ref{asmpVariance1}$, for any $t\in[T]$,

1. if $p_1p_2=o(\min\{Tp_1^2,Tp_2^2\})$,
\begin{equation}\nonumber
\sqrt{p_1p_2}\left[\Vec(\widehat{\Fb}_t)-\left(\widehat{\Hb}_{\Cb}^{-1}\otimes\widehat{\Hb}_{\Rb}^{-1}\right)\Vec(\Fb_t)\right]\stackrel{d}\rightarrow \mathcal{N}(\textbf{0},(\bGamma_2\otimes\bGamma_1)^\top\bSigma^\ast(\bGamma_2\otimes\bGamma_1));
\end{equation}

2. if $p_1p_2\gtrsim\min\{Tp_1^2,Tp_2^2\}$,
\begin{equation}\nonumber
\Vec(\widehat{\Fb}_t)-\left(\widehat{\Hb}_{\Cb}^{-1}\otimes\widehat{\Hb}_{\Rb}^{-1}\right)\Vec(\Fb_t)=O_p\left(\frac{1}{Tp_1^2}+\frac{1}{Tp_2^2}\right).
\end{equation}
\end{theorem}
\begin{theorem}\label{theo3.3}
Under Assumptions $\ref{asmpa-mixing}$ to $\ref{asmpVariance1}$, assume that $\displaystyle\frac{\delta_{pT}}{\sqrt{Tp_1}}$, $\displaystyle\frac{\delta_{pT}}{\sqrt{Tp_2}}$ and $\displaystyle\frac{\delta_{pT}}{\sqrt{p_1p_2}}$ converge to $a$, $b$ and $c$, respectively, where $\delta_{pT}=\min\{\sqrt{Tp_1},\sqrt{Tp_2},\sqrt{p_1p_2}\}$. Then as $\min\{T,p_1,p_2\}\rightarrow\infty$, for any $t\in[T]$, $i\in[p_1]$ and $j\in[p_2]$,  we have
\begin{equation}\nonumber
\frac{\widehat{\textbf{S}}_{t,ij}-\textbf{S}_{t,ij}}{\displaystyle\sqrt{\frac{1}{Tp_1}Z_{t,ij}^{(1)}+\frac{1}{Tp_2}Z_{t,ij}^{(2)}+\frac{1}{p_1p_2}Y_{ij}}}\stackrel{d}\rightarrow N(0,1),
\end{equation}
where $Z_{t,ij}^{(1)}=\Rb_{i\cdot}^\top\Fb_t\bSigma_2^{-1}\Vb_{2j}\bSigma_2^{-1}\Fb_t^\top\Rb_{i\cdot}$, $Z_{t,ij}^{(2)}=\Cb_{j\cdot}^\top\Fb_t^\top\bSigma_1^{-1}\Vb_{1i}\bSigma_1^{-1}\Fb_t\Cb_{j\cdot}$, $Y_{ij}=(\Cb_{j\cdot}\otimes\Rb_{i\cdot})^\top(\bGamma_2\otimes\bGamma_1)^\top\bSigma^\ast(\bGamma_2\otimes\bGamma_1)(\Cb_{j\cdot}\otimes\Rb_{i\cdot})$, and $\widehat{\textbf{S}}_{t,ij}$ and $\textbf{S}_{t,ij}$ are the $(i,j)$th elements of $\widehat{\textbf{S}}$ and $\textbf{S}$, respectively.
\end{theorem}

Theorem $\ref{theo3.2}$ shows that $\vec(\widehat{\Fb}_t)$ is a consistent estimator of $\vec(\Fb_t)$ up to a rotation matrix, and the asymptotic normality holds if $p_1p_2=o(\min\{Tp_1^2,Tp_2^2\})$ and the convergence rates are the same with those of PE estimators derived in \cite{yu2022projected}. Theorem $\ref{theo3.3}$ gives the  asymptotic normality of the one-step estimator for the element-wise common components.

\section{Data-driven Generalized Principal Component Analysis }\label{Section3}
In this section, we introduce our compeletely data-driven GPCA method when the covariance matrices are unknown and need to be estimated. 
In Section \ref{Section2}, we assume that $\Ub$ and $\Vb$ are both known symmetric matrices. In practice, however, both $\Ub$ and $\Vb$ are unknown and need be estimated before implementing the Oracle GPCA procedure in Algorithm \ref{algori2.1}. In this study, we propose  adaptive thresholding estimators for $\Ub$ and $\Vb$, by assuming that $\Ub$ and $\Vb$ are approximately sparse as precisely stated in the following Assumption \ref{SparsityUV}, see similar assumptions adopted also in \cite{bickel2008covariance}, \cite{rothman2009generalized} and \cite{fan2013large} for estimating high-dimensional covariance matrices.

\begin{asmp}[Sparsity of $\Ub$ and $\Vb$]\label{SparsityUV}
There exist $q_r\in[0,1)$ and $q_c\in[0,1)$, such that
$$
m_{p_1}=\max_{i\leq p_1}\sum_{j\leq p_1}|u_{ij}|^{q_r}\quad\textrm{and}\quad m_{p_2}=\max_{i\leq p_2}\sum_{j\leq p_2}|v_{ij}|^{q_c},
$$
do not grow too fast as $p_1\rightarrow\infty$ and $p_2\rightarrow\infty$, where $u_{ij}$ and $v_{ij}$ are the $(i,j)$th elements of $\Ub$ and $\Vb$, respectively.
\end{asmp}

In the following, we introduce in detail how to estimate the covariance matrices $\Ub$ and $\Vb$. Firstly,
we use the PE method by \cite{yu2022projected} to obtain the estimators of loading matrices, denoted by $\ddot{\Rb}$ and $\ddot{\Cb}$. The estimators of the factor matrix and common component  are given by $\ddot{\Fb}_t=\ddot{\Rb}^\top\Xb_t\ddot{\Cb}/(p_1p_2)$ and $\ddot{\mathbf{S}}_t=\ddot{\Rb}\Fb_t\ddot{\Cb}^\top$, respectively. Then, the estimator of idiosyncratic component is  $\ddot{\Eb}_t=\Xb_t-\ddot{\mathbf{S}}_t$ and the sample row and column covariance matrices of the residual matrix time series are defined as
$$
\widetilde{\Ub}=\frac{1}{Tp_2}\sum_{t=1}^T\ddot{\Eb}_t\ddot{\Eb}_t^\top,\quad\textrm{and}\quad\widetilde{\Vb}=\frac{1}{Tp_1}\sum_{t=1}^{T}\ddot{\Eb}_t^\top\ddot{\Eb}_t,
$$
respectively.
Let $\widetilde{u}_{ij}$ be the $(i,j)$th element of $\widetilde{\Ub}$ and  define
$$
\widetilde{\theta}_{ij}^r=\frac{1}{Tp_2}\sum_{t=1}^T\sum_{l=1}^{p_2}\left(\ddot{e}_{t,il}\ddot{e}_{t,jl}-\widetilde{u}_{ij}\right)^2, \quad \displaystyle\omega_T^r=\sqrt{\frac{\log{p_1}}{Tp_2}}+\frac{1}{\sqrt{Tp_1}}+\frac{1}{\sqrt{p_1p_2}}.
$$
The soft thresholding estimator for $\Ub$ is defined as follows:
$$
\widehat{\Ub}=\left(\widehat{u}_{ij}\right)_{p_1\times p_1},\quad\widehat{u}_{ij}=\left\{
\begin{aligned}
\widetilde{u}_{ii},\quad\quad\quad i=j\\
s_{ij}^r\left(\widetilde{u}_{ij}\right),\quad i\neq j
\end{aligned}
\right.,
$$
where $s_{ij}^r(z)=\textrm{sgn}(z)\left(|z|-\tau_{ij}^r\right)_+$, $\tau_{ij}^r=C_r\sqrt{\widetilde{\theta}_{ij}^r}\omega_T^r$, with $C_r>0$ being a large enough constant. 

In parallel, let $\widetilde{v}_{ij}$ be the $(i,j)$th element of $\widetilde{\Vb}$ and define
$$
\widetilde{\theta}_{ij}^c=\frac{1}{Tp_1}\sum_{t=1}^T\sum_{l=1}^{p_1}\left(\ddot{e}_{t,li}\ddot{e}_{t,lj}-\widetilde{v}_{ij}\right)^2,\quad\omega_T^c=\sqrt{\frac{\log{p_2}}{Tp_1}}+\frac{1}{\sqrt{Tp_2}}+\frac{1}{\sqrt{p_1p_2}},
$$
the soft thresholding estimator for $\Vb$ is
$$
\widehat{\Vb}=\left(\widehat{v}_{ij}\right)_{p_2\times p_2},\quad\widehat{v}_{ij}=\left\{
\begin{aligned}
\widetilde{v}_{ii},\quad\quad\quad i=j\\
s_{ij}^r\left(\widetilde{v}_{ij}\right),\quad i\neq j
\end{aligned}
\right.,
$$
where $\tau_{ij}^c=C_c\sqrt{\widetilde{\theta}_{ij}^c}\omega_T^c$ and $s_{ij}^c(z)=\textrm{sgn}(z)\left(|z|-\tau_{ij}^c\right)_+$.

Inspired by \cite{fan2013large}, $C_r$ and $C_c$ can be determined in a data-driven manner and chosen through multifold cross-validation. Firstly, randomly divide $\{\ddot{\Eb}_t\}_{t=1}^T$ into two subsets,  denoted as $\{\ddot{\Eb}_t\}_{t\in J_1}$ and $\{\ddot{\Eb}_t\}_{t\in J_2}$. The sizes of $J_1$ and $J_2$, denoted by $T(J_1)$ and $T(J_2)$, satisfy $T(J_1)\asymp T$ and $T(J_1)+T(J_2)=T$. \cite{bickel2008covariance} suggested to choose $T(J_1)=T(1-(\textrm{log} T)^{-1})$ in sparse matrix estimation. Secondly, repeat this procedure $H$ times. At the $h$th split, we denote $\widehat{\Ub}^h_{J_1}(C_r)$ and $\widehat{\Vb}^h_{J_1}(C_c)$ as the soft thresholding estimators with the threshold $\tau_{ij}^r$ and $\tau_{ij}^c$ applied to the training dataset $\{\ddot{\Eb}_t\}_{t\in J_1}$, respectively. We also define $\widehat{\Ub}^h_{J_2}$ and $\widehat{\Vb}^h_{J_2}$ as the sample covariance based on the validation set, where $\displaystyle\widehat{\Ub}^h_{J_2}=\frac{1}{T(J_2)p_2}\sum_{t\in J_2}\ddot{\Eb}_t\ddot{\Eb}_t^\top$ and $\displaystyle\widehat{\Vb}^h_{J_2}=\frac{1}{T(J_2)p_1}\sum_{t\in J_2}\ddot{\Eb}_t^\top\ddot{\Eb}_t$. Then we can choose the constants $C_r^\ast$ and $C_c^\ast$ by minimizing the following cross-validation objective function
\begin{equation}\nonumber
C_r^\ast=\arg\min\limits_{C_{\min}^r\leq C_r\leq M}\frac{1}{H}\sum_{h=1}^H\left\|\widehat{\Ub}^h_{J_1}(C_r)-\widehat{\Ub}^h_{J_2}\right\|_F^2\quad\textrm{and}\quad C_c^\ast=\arg\min\limits_{C_{\min}^c\leq C_c\leq M}\frac{1}{H}\sum_{h=1}^H\left\|\widehat{\Vb}^h_{J_1}(C_c)-\widehat{\Vb}^h_{J_2}\right\|_F^2,
\end{equation}
where $M$ is a large constant such that $\widehat{\Ub}$ and $\widehat{\Vb}$ are diagonal, $C_{\min}^r$ and $C_{\min}^c$ are the minimum constants that guarantee the positive definiteness of $\widehat{\Ub}$ and $\widehat{\Vb}$, respectively. The resulting $C_r$ and $C_c$ are data-driven, relying on $\Xb_t$ as well as $p_1$, $p_2$ and $T$.  On the other hand, given data matrix $\Xb_t\in\RR^{p_1\times p_2}$, $C_r^\ast$ and $C_c^\ast$ are universal constants in the thresholds $\tau_{ij}^r=C_r^\ast\sqrt{\widetilde{\theta}_{ij}^r}\omega_T^r$ and $\tau_{ij}^c=C_c\sqrt{\widetilde{\theta}_{ij}^c}\omega_T^c$ in the sense that they do not change with respect to the position $(i,j)$. As pointed out in \cite{fan2013large}, cross-validation relies on the covariance matrix of the vectorized idiosyncratic error matrix, rather than that of the vectorized observation matrix. Thus this cross-validation procedure improves the selection of thresholds.

After obtaining the estimators of the covariance matrices of the idiosyncratic component, i.e., $\widehat{\Ub}$ and $\widehat{\Vb}$, we replace $\Ub$ and $\Vb$ for Oracle GPCA method in Section \ref{Section2} with $\hat{\Ub}$ and $\hat{\Vb}$, and perform the same procedure as in Algorithm \ref{algori2.1} to obtain the data-driven GPCA estimators. In detail, we first use $\alpha$-PCA method by \cite{chen2023statistical} with $\alpha=0$ to obtain the initial estimators of the loading matrices $\widetilde{\Rb}_0^\ast$ and $\widetilde{\Cb}_0^\ast$. Let $\widetilde{\Rb}_0=\widehat{\Ub}^{1/2}\widetilde{\Rb}_0^\ast$, $\widetilde{\Cb}_0=\widehat{\Vb}^{1/2}\widetilde{\Cb}_0^\ast$ and define
$$
\widetilde{\Mb}_{\Cb}=\frac{1}{Tp_2}\sum_{t=1}^T\widehat{\Ub}^{-1/2}\Xb_t\widehat{\Vb}^{-1}\widetilde{\Cb}_0\widetilde{\Cb}_0^\top\widehat{\Vb}^{-1}\Xb_t^\top\widehat{\Ub}^{-1/2},$$
and $$\quad\widetilde{\Mb}_{\Rb}=\frac{1}{Tp_1}\sum_{t=1}^T\widehat{\Vb}^{-1/2}\Xb_t^\top\widehat{\Ub}^{-1}\widetilde{\Rb}_0\widetilde{\Rb}_0^\top\widehat{\Ub}^{-1}\Xb_t\widehat{\Vb}^{-1/2}.
$$

We denote the first $k_1$ eigenvectors of $\widetilde{\Mb}_{\Cb}$ as $\{\tilde{\br}(1),\ldots,\tilde{\br}(k_1)\}$ and the corresponding eigenvalues as $\{\tilde{\theta}_1^,\ldots,\tilde{\theta}_{k_1}\}$. Similarly, denote the first $k_2$ eigenvectors of $\widetilde{\Mb}_{\Rb}$ as $\{\tilde{\bc}(1),\ldots,\tilde{\bc}(k_2)\}$ and the corresponding eigenvalues as $\{\tilde{\lambda}_1,\ldots,\tilde{\lambda}_{k_2}\}$. Similarly as in Section $\ref{Section2}$, we have $\widetilde{\Rb}=\sqrt{p_1}\widehat{\Ub}^{1/2}(\tilde{\br}(1),\ldots,\tilde{\br}(k_1))$ and $\widetilde{\Cb}=\sqrt{p_2}\widehat{\Vb}^{1/2}(\tilde{\bc}(1),\ldots,\tilde{\bc}(k_2))$. However, $\widetilde{\Mb}_{\Cb}$ still depends on the unknown column factor loading $\Cb$ while $\widetilde{\Mb}_{\Rb}$ depends on the unknown row factor loading $\Rb$. We also consider an iterative procedure to get the data-driven GPCA estimators and summarize the data-driven GPCA method in Algorithm $\ref{algori3.1}$.
\begin{algorithm}
\setstretch{1.5}
  \caption{Data-driven Generalized Principal Component Analysis }
  \label{algori3.1}
  \begin{algorithmic}[1]
    \Require Data matrices $\{\Xb_t\}_{1\leq t\leq T}$, the estimated covariance matrices $\widehat{\Ub}$ and $\widehat{\Vb}$, the pair of factor numbers $(k_1,k_2)$
    \Ensure Factor loading matrices $\widetilde{\Rb}$ and $\widetilde{\Cb}$
    \State define $\widetilde{\Zb}_t=\widehat{\Ub}^{-1/2}\Xb_t\widehat{\Vb}^{-1/2}$ and obtain the initial estimators $\widetilde{\Rb}_0^\ast$ and $\widetilde{\Cb}_0^\ast$ by $\alpha$-PCA with $\alpha=0$;
    \State define $\widetilde{\Mb}_{\Cb}=(Tp_2)^{-1}\sum_{t=1}^T\widetilde{\Zb}_t\widetilde{\Cb}_0^\ast\widetilde{\Cb}_0^{\ast\top}\widetilde{\Zb}_t^\top$ and $\widetilde{\Mb}_{\Rb}=(Tp_1)^{-1}\sum_{t=1}^T\widetilde{\Zb}_t^\top\widetilde{\Rb}_0^\ast\widetilde{\Rb}_0^{\ast\top}\widetilde{\Zb}_t$, obtain the leading $k_1$ eigenvectors of $\widetilde{\Mb}_{\Cb}$  and the leading $k_2$ eigenvectors of $\widetilde{\Mb}_{\Rb}$, denoted as $\widetilde{\Qb}_R$ and $\widetilde{\Qb}_C$, respectively;
    \State the row and column loading matrices are finally given by
    $\widetilde{\Rb}=\sqrt{p_1}\hat{\Ub}^{1/2}\widetilde{\Qb}_\Rb$ and $\widetilde{\Cb}=\sqrt{p_2}\hat{\Vb}^{1/2}\widetilde{\Qb}_\Cb$.
  \end{algorithmic}
\end{algorithm}

Similarly, as long as the loading matrices are determined, the factor matrix $\Fb_t$ can be estimated by $\widetilde{\Fb}_t=\widetilde{\Rb}^\top\widehat{\Ub}^{-1}\Xb_t\widehat{\Vb}^{-1}\widetilde{\Cb}/(p_1p_2)$. The estimator of the common component is then given by $\widetilde{\textbf{S}}_t=\widetilde{\Rb}\widetilde{\Fb}_t\widetilde{\Cb}^\top$. To establish the theoretical properties of $\widehat{\Ub}$ and $\widehat{\Vb}$, the following assumptions are required.

\begin{asmp}\label{exasum}\ 

1. The elements of $\Eb_t^\ast$ are independent and identically distributed sub-Gaussian random variables with zero mean and unit variance.

2. (1). for any $t\in[T],i\in[p_1],j\in[p_2]$,
\begin{equation}\nonumber
\sum_{s=1}^T\sum_{l=1}^{p_1}\sum_{h=1}^{p_2}\vert\EE \left(e_{t,ij}^\ast e_{s,lh}^\ast\right)\vert\leq c\ .
\end{equation}

(2). for any $t\in[T],i,l_1\in[p_1],j,h_1\in[p_2]$,
\begin{equation}\nonumber
\begin{aligned}
&(i). \sum_{s=1}^T\sum_{l_2=1}^{p_1}\sum_{h=1}^{p_2}\Big{\vert}\textrm{Cov}(e_{t,ij}^\ast e_{t,l_1j}^\ast,e_{s,ih}^\ast e_{s,l_2h}^\ast)\Big{\vert}\leq c,\quad\sum_{s=1}^T\sum_{l=1}^{p_1}\sum_{h_2=1}^{p_2}\Big{\vert}\textrm{Cov}(e_{t,ij}^\ast e_{t,ih_1}^\ast,e_{s,lj}^\ast e_{s,lh_2}^\ast)\Big{\vert}\leq c\ ;\\
&(ii). \sum_{s=1}^T\sum_{l_2=1}^{p_1}\sum_{h_2=1}^{p_2}\left(\Big{\vert}\textrm{Cov}(e_{t,ij}^\ast e_{t,l_1h_1}^\ast,e_{s,ij}^\ast e_{s,l_2h_2}^\ast)\Big{\vert}+\Big{\vert}\textrm{Cov}(e_{t,l_1j}^\ast e_{t,ih_1}^\ast,e_{s,l_2j}^\ast e_{s,ih_2}^\ast)\Big{\vert}\right)\leq c \ .
\end{aligned}
\end{equation}

3. There exists a constant $r_2>0$ with $r_1^{-1}+3r_2^{-1}>1$, and $b_1>0$ such that for any $s>0$, $i\in[p_1]$ and $j\in[p_2]$,
$$
P(\vert f_{t,ij}\vert>s)\leq\exp(-(s/b_1)^{r_2}),
$$
where $f_{t,ij}$ is the $(i,j)$th element of $\Fb_t$.

4. For all $t\in[T]$, $i\in[p_1]$, $j\in[p_2]$, $k\in[k_1]$ and $l\in[k_2]$, $\EE \left(e^\ast_{t,ij}f_{t,kl}\right)$=0.
\end{asmp}

{
Similar to Assumption $\ref{asmpNoise}$, Assumption $\ref{exasum}$ 2 still allows weak temporal correlations of the noises. Meanwhile, Assumptions $\ref{exasum}$ 1 and 3 impose exponential-type tails on factors and noises, which enable the application of large deviation theory to establish the theoretical properties of $\widehat{\Ub}$ and $\widehat{\Vb}$. Although Assumption $\ref{exasum}$ 4 is slightly stronger than Assumption $\ref{asmpFactorNoise}$, it remains standard in the high-dimensional covariance matrix estimation literature and helps simplify the technical analysis. Importantly, the theoretical results of \cite{yu2022projected} still hold, even with the addition of Assumption $\ref{exasum}$. The following theorem establishes the convergence rates of the estimators of idiosyncratic covariance matrices.
}
\begin{theorem}\label{theoUV}
Under Assumptions $\ref{asmpa-mixing}$ to $\ref{asmpLoading}$, Assumption $\ref{asmpUV}$ and Assumptions $\ref{SparsityUV}$ to $\ref{exasum}$, we have
$$
\left\|\widehat{\Ub}-\frac{\tr(\Vb)}{p_2}\Ub\right\|=O_p\left((\omega_T^r)^{1-q_r}m_{p_1}\right)\quad\textrm{and}\quad\left\|\widehat{\Vb}-\frac{\tr(\Ub)}{p_1}\Vb\right\|=O_p\left((\omega_T^c)^{1-q_c}m_{p_2}\right).
$$
If further $(\omega_T^r)m_{p_1}=o(1)$ and $(\omega_T^c)m_{p_2}=o(1)$, then $\widehat{\Ub}$ and $\widehat{\Vb}$ are invertible with probability approaching one, and
$$
\left\|(\widehat{\Ub})^{-1}-\left(\frac{\tr(\Vb)}{p_2}\Ub\right)^{-1}\right\|=O_p\left((\omega_T^r)^{1-q_r}m_{p_1}\right)\quad\textrm{and}\quad\left\|(\widehat{\Vb})^{-1}-\left(\frac{\tr(\Ub)}{p_1}\Vb\right)^{-1}\right\|=O_p\left((\omega_T^c)^{1-q_c}m_{p_2}\right).
$$
\end{theorem}

{Theorem $\ref{theoUV}$ is a  matrix form generalization of the Theorem 3.1 in \cite{fan2013large}, 
indicating that $\widehat{\Ub}$ and $\widehat{\Vb}$ are the consistent estimators of $\displaystyle{\tr(\Vb)}\Ub/{p_2}$ and $\displaystyle{\tr(\Ub)}\Vb/{p_1}$, respectively. The covariance matrices $\Ub$ and $\Vb$ are not separately identifiable. Indeed,  $\widehat{\Ub}$ ($\widehat{\Vb}$) is an estimator of the true covariance matrix $\Ub$ ($\Vb$) up to a constant. However, this would not bring any trouble for the estimation of the loading spaces and the common component  by our GPCA method. It can be proved that the estimated covariance matrix of the vectorized idiosyncratic error matrix is $\displaystyle{(Tp_1p_2\widehat{\Vb}\otimes\widehat{\Ub})}/{(\sum_{t=1}^T\|\ddot{\Eb}_t\|_F^2})$. Thus, we can assign constant $\displaystyle{(Tp_1p_2)}/{(\sum_{t=1}^T\|\ddot{\Eb}_t\|_F^2})$ to either $\widehat{\Ub}$ or $\widehat{\Vb}$ for identifiability issue.
The following theorem presents the asymptotic distributions of the data-driven GPCA estimators.}
\begin{theorem}\label{theo3.4}
Under Assumptions $\ref{asmpa-mixing}$ to $\ref{asmpLoading}$, Assumptions $\ref{asmpUV}$ to $\ref{asmpVariance1}$ and Assumptions $\ref{SparsityUV}$ to $\ref{exasum}$, if the initial estimators $\widetilde{\Rb}_0^{\ast}$ and $\widehat{\Cb}_0^{\ast}$ are obtained by $\alpha$-PCA as stated in Algorithm $\ref{algori3.1}$, there exist matrices $\widetilde{\Hb}_{\Rb}$ and $\widetilde{\Hb}_{\Cb}$, satisfying $\widetilde{\Hb}_{\Rb}^\top\widetilde{\Hb}_{\Rb}\stackrel{p}\rightarrow\Ib_{k_1}$ and $\widetilde{\Hb}_{\Cb}^\top\widetilde{\Hb}_{\Cb}\stackrel{p}\rightarrow\Ib_{k_2}$, such that

1. for $i\in[p_1]$,
\begin{equation}\nonumber
\begin{cases}
\sqrt{Tp_2}(\tilde{\Rb}_{i\cdot}-\tilde{\Hb}_{\Rb}^\top\Rb_{i\cdot})\stackrel{d}\rightarrow\mathcal{N}(\textbf{0},\bLambda_1^{-1}\bGamma_1^\top\Ab_{1i}\bGamma_1\bLambda_1^{-1}),\quad\textrm{if}\quad Tp_2=o(\min\{T^2p_1^2,p_1^2p_2^2\}),\\
\displaystyle\tilde{\Rb}_{i\cdot}-\tilde{\Hb}_{\Rb}^\top\Rb_{i\cdot}=O_p\left(\frac{1}{Tp_1}+\frac{1}{p_1p_2}\right),\quad\quad\quad\quad\quad\quad\quad\textrm{if}\quad Tp_2\gtrsim\min\{T^2p_1^2,p_1^2p_2^2\},
\end{cases}
\end{equation}

 for $j\in[p_2]$,
\begin{equation}\nonumber
\begin{cases}
\sqrt{Tp_1}(\tilde{\Cb}_{j\cdot}-\tilde{\Hb}_{\Cb}^\top\Cb_{j\cdot})\stackrel{d}\rightarrow\mathcal{N}(\textbf{0},\bLambda_2^{-1}\bGamma_2^\top\Ab_{2j}\bGamma_2\bLambda_2^{-1}),\quad\textrm{if}\quad Tp_1=o(\min\{T^2p_2^2,p_1^2p_2^2\}),\\
\displaystyle\tilde{\Cb}_{j\cdot}-\tilde{\Hb}_{\Cb}^\top\Cb_{j\cdot}=O_p\left(\frac{1}{Tp_2}+\frac{1}{p_1p_2}\right),\quad\quad\quad\quad\quad\quad\quad\textrm{if}\quad Tp_1\gtrsim\min\{T^2p_2^2,p_1^2p_2^2\},
\end{cases}
\end{equation}
where $\widetilde{\Rb}_{i\cdot}$ and $\widetilde{\Cb}_{j\cdot}$ are the $i$th/$j$th row of $\widetilde{\Rb}$ and $\widetilde{\Cb}$,  $\widetilde{\Rb}$ and $\widetilde{\Cb}$ are the one-step estimators of $\Rb$ and $\Cb$ respectively;

2. for $t\in[T]$,
\begin{equation}\nonumber
\begin{cases}
\sqrt{p_1p_2}\left(\Vec(\tilde{\Fb}_t)-\left(\tilde{\Hb}_{\Cb}^{-1}\otimes\tilde{\Hb}_{\Rb}^{-1}\right)\Vec(\Fb_t)\right)\stackrel{d}\rightarrow \mathcal{N}(\textbf{0},(\bGamma_2\otimes\bGamma_1)^\top\bSigma^\ast(\bGamma_2\otimes\bGamma_1)),\quad\textrm{if}\quad p_1p_2=o(\min\{Tp_1^2,Tp_2^2\}),\\
\displaystyle\Vec(\tilde{\Fb}_t)-\left(\tilde{\Hb}_{\Cb}^{-1}\otimes\tilde{\Hb}_{\Rb}^{-1}\right)\Vec(\Fb_t)=O_p\left(\frac{1}{Tp_1^2}+\frac{1}{Tp_2^2}\right),\quad\quad\quad\ \ \quad \quad \quad\quad\quad\quad\quad\textrm{if}\quad p_1p_2\gtrsim\min\{Tp_1^2,Tp_2^2\};
\end{cases}
\end{equation}

3. for $t\in[T]$, $i\in[p_1]$ and $j\in[p_2]$,
\begin{equation}\nonumber
\frac{\tilde{\textbf{S}}_{t,ij}-\textbf{S}_{t,ij}}{\displaystyle\sqrt{\frac{1}{Tp_1}Z_{t,ij}^{(1)}+\frac{1}{Tp_2}Z_{t,ij}^{(2)}+\frac{1}{p_1p_2}Y_{ij}}}\stackrel{d}\rightarrow N(0,1),
\end{equation}
where $Z_{t,ij}^{(1)}$, $Z_{t,ij}^{(2)}$ and $Y_{ij}$ are defined in Theorem \ref{theo3.3}, $\widetilde{\textbf{S}}_{t,ij}$ is the $(i,j)$th element of $\widetilde{\textbf{S}}$.
\end{theorem}

Theorem $\ref{theo3.4}$ shows that the data-driven GPCA estimators (of the loadings, factors and common components) have the same asymptotic variances as those of the corresponding Oracle GPCA estimators. In other words, Theorem $\ref{theo3.4}$ establishes the equivalence of the asymptotic behavior between Oracle GPCA and data-driven GPCA estimators, though under a slight stronger Assumption $\ref{exasum}$. Thus, discussions after Theorem $\ref{theo3.1}$ on the comparison of estimation efficiency between the Oracle GPCA method and the PE method also apply to the data-driven GPCA method. 
\section{Simulation Study}\label{Section4}
\subsection{Data Generation}
In this section, we introduce the data generation mechanism of the synthetic dataset to verify the finite sample performances of the proposed Oracle GPCA and data-driven GPCA methods in Algorithms $\ref{algori2.1}$ and $\ref{algori3.1}$.

We set $k_1=k_2=3$, and draw the elements of $\Rb$ and $\Cb$ independently from uniform distribution $\mathcal{U}(-1,1)$, and let
$$
\vec(\Fb_t)=\phi\times\vec(\Fb_{t-1})+\sqrt{1-\phi^2}\times\mathbf{\epsilon}_t,\ \mathbf{\epsilon}_t\stackrel{i.i.d.}\sim\mathcal{N}(\textbf{0},\Ib_{k_1\times k_2}),
$$
$$
\vec(\Eb_t)=\psi\times\vec(\Eb_{t-1})+\sqrt{1-\psi^2}\times\vec(\Ub_t),\ \Ub_t\stackrel{i.i.d.}\sim\mathcal{MN}(\textbf{0},\Ub,\Vb),
$$
where $\phi$ and $\psi$ control the temporal correlations. We consider the following three cases to generate $\Ub$ and $\Vb$.

\vspace{0.5em}

$\textbf{Case 1 (Weak Cross-sectionally Correlated Matrix).}$ The diagonal elements of $\Ub$ and $\Vb$ are 1, and the off-diagonal elements are $1/p_1$ and $1/p_2$, respectively.

\vspace{0.5em}

$\textbf{Case 2 (Banded Matrix).}$ The diagonal elements of $\Ub$ and $\Vb$ are 1, and the $i$th row and $j$th column off-diagonal elements of $\Ub$ and $\Vb$ are both $\displaystyle\left(\frac{1}{2}-\frac{|i-j|}{10}\right)_{+}$.

\vspace{0.5em}

$\textbf{Case 3 (Block Dignonal Matrix).}$ Each sub block of $\Ub$ and $\Vb$ is composed of 5 symmetric matrices with dimensions $p_1/5\times p_1/5$ and $p_2/5\times p_2/5$, respectively. The off-diagonal elements of the subblocks are $10/p_1$ and $10/p_2$, and the diagonal elements are randomly selected from the interval [1,1.4] with an interval of 0.1.

\vspace{0.5em}

We compare the Oracle GPCA and data-driven GPCA methods in Algorithms $\ref{algori2.1}$ and $\ref{algori3.1}$, the $\alpha$-PCA method by \cite{chen2023statistical} and the PE method by \cite{yu2022projected}, in terms of estimating the loading spaces and the common component matrices. For convenience and ease of differentiation, we denote ``Oracle" as the Oracle GPCA method and ``GPCA" as the data-driven GPCA method. All the following simulation results are based on 200 repetitions. We set $\phi=\psi=0.1$ and consider the following two settings for the combinations of $(T,p_1,p_2)$:
\begin{itemize}
\item[$\bullet$]\textbf{Setting A:} $p_1=10$, $T=p_2\in\{20,50,100,150,200\}$ .
\item[$\bullet$]\textbf{Setting B:} $p_2=10$, $T=p_1\in\{20,50,100,150,200\}$.
\end{itemize}
\subsection{Estimation Error for Loading Spaces}\label{sec:4.2}
In this section, we compare the performances  of the GPCA, Oracle, $\alpha$-PCA and PE methods in terms of estimating the loading spaces.  For two column-wise orthogonal matrices $\Qb_1\in\RR^{p\times q_1}$ and $\Qb_2\in\RR^{p\times q_2}$, the distance
between the spaces spanned by the columns of  $\Qb_1$ and $\Qb_2$ can be measured by
$$
\mathcal{D}(\Qb_1,\Qb_2)=\left(1-\frac{1}{\min\{q_1,q_2\}}\tr\left(\Qb_1\Qb_1^\top\Qb_2\Qb_2^\top\right)\right)^{1/2}.
$$
By the definition of $\mathcal{D}(\Qb_1,\Qb_2)$, we can easily see that the distance between the column spaces of $\Qb_1$ and $\Qb_2$ is always in the interval $[0,1]$. In particular, when $\Qb_1$ and $\Qb_2$ span the same space, the distance $\mathcal{D}(\Qb_1,\Qb_2)$ is equal to 0, while is equal to 1 when the two spaces are orthogonal. The Gram-Schmidt orthogonalization can be used to make $\Qb_1$ and $\Qb_2$ column-orthogonal matrices if not.
\begin{table}[!ht]
	 	\caption{Averaged estimation errors and standard errors (in parentheses) of $\mathcal{D}(\widetilde{\Rb},\Rb)$ and $\mathcal{D}(\widetilde{\Cb},\Cb)$ for Settings A and B under Case 1 over 200 replications. ``GPCA": data-driven GPCA method. ``Oracle": oracle GPCA method. ``$\alpha$-PCA": the $\alpha$-PCA method with $\alpha=0$. ``PE": the projection estimation method.}
	 	 \label{tab:main1}\renewcommand{\arraystretch}{1} \centering
	 	\selectfont
	 	\begin{threeparttable}
	 		 \scalebox{0.9}{\begin{tabular*}{16.5cm}{ccccccccc}
                    \toprule[2pt]
	 				&\multirow{2}{*}{Evaluation} &\multirow{2}{*}{T}
                    &\multirow{2}{*}{$p_1$}    &\multirow{2}{*}{$p_2$}
                    &\multirow{2}{*}{GPCA}
                    &\multirow{2}{*}{Oracle}
                    &\multirow{2}{*}{$\alpha$-PCA} &\multirow{2}{*}{PE}
                    \cr
	 			     \\
	 				\midrule[1pt]
	 				 &\multicolumn{8}{c}{\multirow{1}{*}{$\textbf{Weak Cross-sectionally Correlated Matrix}$}}\\ &\multirow{5}*{$\mathcal{D}(\widetilde{\Rb},\Rb)$}&20&\multirow{5}*{20}&20&0.1003(0.0218)&0.0918(0.0166)
     &0.2603(0.1362)&0.1021(0.0229)\\
	 				&&50&&50&0.0365(0.0061)&0.0355(0.0058)&0.1971(0.1013)&0.0373(0.0064)\\
                    &&100&&100&0.0183(0.0027)&0.0181(0.0028)&0.1913(0.1045)&0.0187(0.0031)\\
                    &&150&&150&0.0119(0.0020)&0.0118(0.0019)&0.1737(0.0969)&0.0121(0.0021)\\
                    &&200&&200&0.0089(0.0013)&0.0089(0.0013)&0.1790(0.0996)&0.0091(0.0015)\\
                    
                    &\multirow{5}*{$\mathcal{D}(\widetilde{\Cb},\Cb)$}&20&\multirow{5}*{20}&20&0.1003(0.0203)&0.0921(0.0162)&0.2593(0.1160)&0.1023(0.0199)\\
	 				&&50&&50&0.0562(0.0066)&0.0556(0.0064)&0.0716(0.0139)&0.0588(0.0075)\\
                    &&100&&100&0.0385(0.0033)&0.0386(0.0033)&0.0429(0.0047)&0.0407(0.0041)\\
                    &&150&&150&0.0314(0.0025)&0.0315(0.0025)&0.0341(0.0034)&0.0332(0.0031)\\
                    &&200&&200&0.0273(0.0021)&0.0274(0.0021)&0.0293(0.0027)&0.0288(0.0027)\\
                    \hline
                    &\multirow{5}*{$\mathcal{D}(\widetilde{\Rb},\Rb)$}&20&20&\multirow{5}*{20}&0.1003(0.0218)&0.0918(0.0166)&0.2603(0.1362)&0.1021(0.0229)\\
	 	&&50&50&&0.0563(0.0062)&0.0557(0.0058)&0.0705(0.0130)&0.0588(0.0069)\\
                    &&100&100&&0.0385(0.0034)&0.0386(0.0034)&0.0426(0.0046)&0.0404(0.0040)\\
                    &&150&150&&0.0318(0.0025)&0.0320(0.0026)&0.0344(0.0033)&0.0335(0.0031)\\
                    &&200&200&&0.0270(0.0020)&0.0271(0.0020)&0.0289(0.0025)&0.0284(0.0024)\\
                    &\multirow{5}*{$\mathcal{D}(\widetilde{\Cb},\Cb)$}&20&20&\multirow{5}*{20}&0.1003(0.0203)&0.0921(0.0162)&0.2593(0.1160)&0.1023(0.0199)\\
	 				&&50&50&&0.0365(0.0059)&0.0356(0.0056)&0.2095(0.1146)&0.0376(0.0062)\\
                    &&100&100&&0.0181(0.0029)&0.0180(0.0028)&0.1928(0.1110)&0.0185(0.0031)\\
                    &&150&150&&0.0121(0.0018)&0.0121(0.0018)&0.1925(0.1027)&0.0124(0.0020)\\
                    &&200&200&&0.0090(0.0014)&0.0090(0.0014)&0.1799(0.1128)&0.0093(0.0016)\\

	 				\bottomrule[2pt]
 				\end{tabular*}}
 			\end{threeparttable}
     \end{table}

\begin{table}[!ht]
	 	\caption{Averaged estimation errors and standard errors (in parentheses) of $\mathcal{D}(\widetilde{\Rb},\Rb)$ and $\mathcal{D}(\widetilde{\Cb},\Cb)$ for Settings A and B under Case 2 over 200 replications.  ``GPCA": data-driven GPCA method. ``Oracle": oracle GPCA method. ``$\alpha$-PCA": the $\alpha$-PCA method with $\alpha=0$. ``PE": the projection estimation method.}
	 	 \label{tab:main2}\renewcommand{\arraystretch}{1} \centering
	 	\selectfont
	 	\begin{threeparttable}
	 		 \scalebox{0.9}{\begin{tabular*}{16.5cm}{ccccccccc}
                    \toprule[2pt]
	 				&\multirow{2}{*}{Evaluation} &\multirow{2}{*}{T}
                    &\multirow{2}{*}{$p_1$}    &\multirow{2}{*}{$p_2$}
                    &\multirow{2}{*}{GPCA}
                    &\multirow{2}{*}{Oracle}
                    &\multirow{2}{*}{$\alpha$-PCA} &\multirow{2}{*}{PE}
                    \cr
	 			     \\
	 				\midrule[1pt]
	 				 &\multicolumn{8}{c}{\multirow{1}{*}{$\textbf{Banded Matrix}$}}\\
       &\multirow{5}*{$\mathcal{D}(\widetilde{\Rb},\Rb)$}&20&\multirow{5}*{20}&20&0.5466(0.3133)&0.0846(0.0216)&0.7223(0.0690)&0.6673(0.1983)\\
	 				&&50&&50&0.0377(0.0677)&0.0319(0.0068)&0.7248(0.0665)&0.1017(0.0784)\\
                    &&100&&100&0.0155(0.0032)&0.0156(0.0031)&0.7262(0.0586)&0.0410(0.0112)\\
                    &&150&&150&0.0096(0.0019)&0.0102(0.0020)&0.7272(0.0646)&0.0259(0.0061)\\
                    &&200&&200&0.0073(0.0014)&0.0077(0.0014)&0.7381(0.0611)&0.0192(0.0044)\\
                    &\multirow{5}*{$\mathcal{D}(\widetilde{\Cb},\Cb)$}&20&\multirow{5}*{20}&20&0.5484(0.3162)&0.0857(0.0208)&0.7210(0.0678)&0.6660(0.1979)\\
	 				&&50&&50&0.0540(0.0668)&0.0518(0.0091)&0.4763(0.1204)&0.1291(0.0753)\\
                    &&100&&100&0.0334(0.0038)&0.0359(0.0044)&0.1819(0.0336)&0.0698(0.0126)\\
                    &&150&&150&0.0274(0.0028)&0.0297(0.0034)&0.1157(0.0188)&0.0554(0.0095)\\
                    &&200&&200&0.0236(0.0024)&0.0254(0.0028)&0.0863(0.0127)&0.0458(0.0079)\\
                    \hline
                    &\multirow{5}*{$\mathcal{D}(\widetilde{\Rb},\Rb)$}&20&20&\multirow{5}*{20}&0.5466(0.3133)&0.0846(0.0216)&0.7223(0.0690)&0.6673(0.1983)\\
	 				&&50&50&&0.0865(0.1810)&0.0516(0.0092)&0.4853(0.1148)&0.1575(0.1744)\\
                    &&100&100&&0.0333(0.0039)&0.0356(0.0044)&0.1823(0.0335)&0.0683(0.0137)\\
                    &&150&150&&0.0277(0.0029)&0.0297(0.0035)&0.1172(0.0195)&0.0540(0.0107)\\
                    &&200&200&&0.0235(0.0024)&0.0253(0.0028)&0.0852(0.0121)&0.0455(0.0071)\\
                    &\multirow{5}*{$\mathcal{D}(\widetilde{\Cb},\Cb)$}&20&20&\multirow{5}*{20}&0.5484(0.3162)&0.0857(0.0208)&0.7210(0.0678)&0.6660(0.1979)\\
	 				&&50&50&&0.0693(0.1788)&0.0312(0.0065)&0.7268(0.0635)&0.1336(0.1743)\\
                    &&100&100&&0.0155(0.0031)&0.0158(0.0032)&0.7306(0.0633)&0.0400(0.0112)\\
                    &&150&150&&0.0098(0.0021)&0.0103(0.0022)&0.7370(0.0695)&0.0262(0.0062)\\
                    &&200&200&&0.0073(0.0016)&0.0077(0.0017)&0.7345(0.0594)&0.0189(0.0050)\\

	 				\bottomrule[2pt]
 				\end{tabular*}}
 			\end{threeparttable}
     \end{table}

     \begin{table}[!ht]
	 	\caption{Averaged estimation errors and standard errors (in parentheses) of $\mathcal{D}(\widetilde{\Rb},\Rb)$ and $\mathcal{D}(\widetilde{\Cb},\Cb)$ for Settings A and B under Case 3 over 200 replications.  ``GPCA": data-driven GPCA method. ``Oracle": oracle GPCA method. ``$\alpha$-PCA": the $\alpha$-PCA method with $\alpha=0$. ``PE": the projection estimation method.}
	 	 \label{tab:main3}\renewcommand{\arraystretch}{1} \centering
	 	\selectfont
	 	\begin{threeparttable}
	 		 \scalebox{0.9}{\begin{tabular*}{16.5cm}{ccccccccc}
                    \toprule[2pt]
	 				&\multirow{2}{*}{Evaluation} &\multirow{2}{*}{T}
                    &\multirow{2}{*}{$p_1$}    &\multirow{2}{*}{$p_2$}
                    &\multirow{2}{*}{GPCA}
                    &\multirow{2}{*}{Oracle}
                    &\multirow{2}{*}{$\alpha$-PCA} &\multirow{2}{*}{PE}
                    \cr
	 			     \\
	 				\midrule[1pt]
	 				 &\multicolumn{8}{c}{\multirow{1}{*}{$\textbf{Block Dignonal Matrix}$}}\\
       
       &\multirow{5}*{$\mathcal{D}(\widetilde{\Rb},\Rb)$}&20&\multirow{5}*{20}&20&0.8444(0.1690)&0.1200(0.0298)&0.8213(0.0518)&0.8689(0.0692)\\
	 				&&50&&50&0.1732(0.2940)&0.0580(0.0123)&0.7665(0.0673)&0.2833(0.2892)\\
                    &&100&&100&0.0292(0.0058)&0.0297(0.0057)&0.7273(0.0606)&0.0496(0.0115)\\
                    &&150&&150&0.0196(0.0036)&0.0201(0.0036)&0.7255(0.0584)&0.0307(0.0061)\\
                    &&200&&200&0.0147(0.0027)&0.0151(0.0028)&0.7121(0.0629)&0.0223(0.0047)\\
                    &\multirow{5}*{$\mathcal{D}(\widetilde{\Cb},\Cb)$}&20&\multirow{5}*{20}&20&0.8391(0.1706)&0.1170(0.0296)&0.8186(0.0576)&0.8623(0.0692)\\
	 				&&50&&50&0.1612(0.2826)&0.0649(0.0113)&0.7117(0.0956)&0.3142(0.2902)\\
                    &&100&&100&0.0394(0.0046)&0.0429(0.0055)&0.2893(0.0905)&0.0845(0.0150)\\
                    &&150&&150&0.0309(0.0035)&0.0339(0.0044)&0.1349(0.0310)&0.0634(0.0106)\\
                    &&200&&200&0.0266(0.0027)&0.0292(0.0034)&0.0887(0.0166)&0.0535(0.0076)\\
                    \hline
                    &\multirow{5}*{$\mathcal{D}(\widetilde{\Rb},\Rb)$}&20&20&\multirow{5}*{20}&0.8444(0.1690)&0.1200(0.0298)&0.8213(0.0518)&0.8689(0.0692)\\
	 				&&50&50&&0.2468(0.3656)&0.0646(0.0919)&0.7228(0.0919)&0.3908(0.3393)\\
                    &&100&100&&0.0435(0.0676)&0.0424(0.0054)&0.2900(0.0881)&0.0892(0.0655)\\
                    &&150&150&&0.0312(0.0035)&0.0341(0.0042)&0.1378(0.0293)&0.0636(0.0100)\\
                    &&200&200&&0.0267(0.0028)&0.0293(0.0035)&0.0880(0.0156)&0.0537(0.0079)\\
                    &\multirow{5}*{$\mathcal{D}(\widetilde{\Cb},\Cb)$}&20&20&\multirow{5}*{20}&0.8391(0.1706)&0.1170(0.0296)&0.8186(0.0576)&0.8623(0.0692)\\
	 				&&50&50&&0.2369(0.3520)&0.0574(0.0112)&0.7605(0.0658)&0.3515(0.3362)\\
                    &&100&100&&0.0340(0.0658)&0.0298(0.0056)&0.7267(0.0632)&0.0541(0.0646)\\
                    &&150&150&&0.0194(0.0036)&0.0199(0.0037)&0.7209(0.0632)&0.0308(0.0067)\\
                    &&200&200&&0.0148(0.0026)&0.0151(0.0028)&0.7100(0.0722)&0.0222(0.0048)\\

	 				\bottomrule[2pt]
 				\end{tabular*}}
 			\end{threeparttable}
     \end{table}
Tables $\ref{tab:main1}$, $\ref{tab:main2}$ and $\ref{tab:main3}$ show the averaged estimation errors with standard errors in parentheses under Case 1, 2 and 3, respectively. First we focus on the case that $\Ub$ and $\Vb$ are weak cross-sectionally correlated matrices. All methods benefit from large dimensions. The GPCA, Oracle and PE methods perform similarly and consistently show significant  advantages over $\alpha$-PCA in terms of estimating $\Rb$ and $\Cb$ when $p_1$ and $p_2$ are small, which is consistent with the conclusion of \cite{yu2022projected}. Then we move to the case where $\Ub$ and $\Vb$ are banded matrices. Obviously, the GPCA and Oracle methods show great advantage over  $\alpha$-PCA and PE methods. The estimation errors by PE method is at least twice of those by GPCA and Oracle methods, which indicates that GPCA and Oracle methods fully utilize the information of sample covariances to improve the estimation accuracy. When $\Ub$ and $\Vb$ are block diagonal matrices,  GPCA and Oracle methods outperform  both  $\alpha$-PCA method and PE method, which further demonstrates the superiority of GPCA and Oracle methods when the idiosyncratic covariance matrices have more complex structures. In addition, the Oracle and the GPCA methods always perform comparably, which indicates that estimating the covariance matrices hardly bring any loss for estimating the loading spaces.

\subsection{Estimation Error for Common Components}
In this section, we will compare the performances of the GPCA method  and Oracle method with those of the $\alpha$-PCA and PE methods in terms of estimating the common component matrices. We evaluate the performances of different methods by the Mean Squared Error, i.e.,
$$
\textrm{MSE}=\frac{1}{Tp_1p_2}\sum_{t=1}^T\left\|\widetilde{\textbf{S}}_t-\textbf{S}_t\right\|_F^2.
$$
     \begin{table}[!ht]
	 	\caption{Averaged mean squared errors and standard errors (in parentheses) for Settings A and B under Case1, 2 and 3 over 200 replications.  ``GPCA": data-driven GPCA method. ``Oracle": oracle GPCA method. ``$\alpha$-PCA": the $\alpha$-PCA method with $\alpha=0$. ``PE": the projection estimation method.}
	 	 \label{tab:main4}\renewcommand{\arraystretch}{1} \centering
	 	\selectfont
	 	\begin{threeparttable}
	 		 \scalebox{0.9}{\begin{tabular*}{13.5cm}{ccccccccc}
                    \toprule[2pt]
	 				&\multirow{2}{*}{Case} &\multirow{2}{*}{T}
                    &\multirow{2}{*}{GPCA}
                    &\multirow{2}{*}{Oracle}
                    &\multirow{2}{*}{$\alpha$-PCA} &\multirow{2}{*}{PE}
                    \cr
	 			     \\
	 				\midrule[1pt]
	 				 &\multicolumn{6}{c}{\multirow{1}{*}{Setting A: $p_1=20$, $T=p_2$}}\\
       
       &\multirow{5}*{1}&20&0.0424(0.0078)&0.0357(0.0040)&1.0036(0.2142)&0.0444(0.0080)\\
	 				&&50&0.0132(0.0013)&0.0126(0.0011)&1.0005(0.1691)&0.0144(0.0019)\\
                    &&100&0.0061(0.0004)&0.0060(0.0004)&1.0011(0.1287)&0.0067(0.0008)\\
                    &&150&0.0039(0.0002)&0.0039(0.0002)&0.9991(0.1305)&0.0044(0.0005)\\
                    &&200&0.0029(0.0001)&0.0029(0.0002)&0.9924(0.1222)&0.0032(0.0003)\\
                    &\multirow{5}*{2}&20&0.9771(0.7466)&0.0222(0.0061)&1.0064(0.2145)&1.2663(0.5103)\\
	 				&&50&0.0128(0.0957)&0.0065(0.0013)&1.0015(0.1691)&0.0652(0.1037)\\
                    &&100&0.0024(0.0004)&0.0029(0.0005)&1.0015(0.1287)&0.0195(0.0058)\\
                    &&150&0.0015(0.0002)&0.0019(0.0003)&0.9994(0.1305)&0.0119(0.0031)\\
                    &&200&0.0011(0.0002)&0.0014(0.0002)&0.9927(0.1222)&0.0084(0.0024)\\
                    &\multirow{5}*{3}&20&2.4168(0.6358)&0.0442(0.0155)&0.9965(0.1914)&2.5740(0.2919)\\
	 				&&50&0.2037(0.5016)&0.0157(0.0035)&1.0064(0.1583)&0.3765(0.5457)\\
                    &&100&0.0062(0.0010)&0.0073(0.0014)&1.0125(0.1374)&0.0284(0.0073)\\
                    &&150&0.0038(0.0006)&0.0046(0.0010)&0.9983(0.1243)&0.0160(0.0040)\\
                    &&200&0.0029(0.0005)&0.0035(0.0008)&1.0005(0.1245)&0.0116(0.0027)\\
                    \hline
                    &\multicolumn{6}{c}{\multirow{1}{*}{Setting B: $p_2=20$, $T=p_1$}}\\
                  &\multirow{5}*{1}&20&0.0424(0.0078)&0.0357(0.0040)&1.0036(0.2142)&0.0444(0.0080)\\
	 		&&50&0.0132(0.0014)&0.0126(0.0010)&1.0041(0.1610)&0.0144(0.0018)\\
                    &&100&0.0060(0.0004)&0.0060(0.0004)&0.1000(0.1395)&0.0066(0.0007)\\
                    &&150&0.0039(0.0002)&0.0040(0.0002)&0.9822(0.1307)&0.0044(0.0005)\\
                    &&200&0.0029(0.0001)&0.0029(0.0001)&1.0019(0.1251)&0.0032(0.0003)\\
                    &\multirow{5}*{2}&20&0.9771(0.7466)&0.0222(0.0061)&1.0064(0.2145)&1.2663(0.5103)\\
	 				&&50&0.0589(0.2602)&0.0064(0.0014)&1.0051(0.1610)&0.1136(0.2655)\\
                    &&100&0.0024(0.0003)&0.0028(0.0005)&1.0004(0.1395)&0.0188(0.0060)\\
                    &&150&0.0015(0.0002)&0.0019(0.0003)&0.9824(0.1307)&0.0116(0.0036)\\
                    &&200&0.0011(0.0001)&0.0014(0.0002)&1.0021(0.1251)&0.0084(0.0023)\\
                    &\multirow{5}*{3}&20&2.4168(0.6358)&0.0442(0.0155)&0.9965(0.1914)&2.5740(0.2919)\\
	 				&&50&0.3389(0.6408)&0.0158(0.0034)&0.9979(0.1597)&0.5301(0.6543)\\
                    &&100&0.0125(0.0895)&0.0072(0.0014)&1.0087(0.1412)&0.0349(0.0892)\\
                    &&150&0.0039(0.0006)&0.0047(0.0010)&1.0002(0.1241)&0.0162(0.0040)\\
                    &&200&0.0029(0.0005)&0.0035(0.0007)&1.0010(0.1207)&0.0116(0.0029)\\
	 							
	 				\bottomrule[2pt]
 				\end{tabular*}}
 			\end{threeparttable}
     \end{table}
Table $\ref{tab:main4}$ shows the averaged MSEs with standard errors in parentheses under Setting A and B for Case 1, 2 and 3. From Table $\ref{tab:main4}$, we see that the GPCA, Oracle and PE methods perform comparably in Case 1, and  outperform the $\alpha$-PCA method. In contrast, when $\Ub$ and $\Vb$ are banded matrices or block dignonal matrices, the GPCA method and Oracle method exhibit better performances than the PE method and $\alpha$-PCA method in terms of estimating the common components. Note that the GPCA method and Oracle method are comparable in all cases, which is consistent with the conclusions drawn for the factor loadings in  Section $\ref{sec:4.2}$.
\subsection{Verifying the Asymptotic Normality}
In this section, we verify the asymptotic normality of $\widetilde{\Rb}$ and  the asymptotic variances obtained in Theorems $\ref{theo3.1}$ and $\ref{theo3.4}$ by numerical studies. For data generation, we set $k_1=k_2=3$, $p_1=20$ and $T=p_2=200$. We consider Case 1 and 2 to generate $\Ub$ and $\Vb$ and draw the elements of $\Rb$ and $\Cb$ independently from uniform distribution $\mathcal{U}(-1,1)$. In addition, we normalize $\Rb^\ast$ as $\sqrt{p_1}$ times its left singular-vector matrix such that the identification condition $\Rb^{\ast\top}\Rb^\ast/p_1=\Ib_{k_1}$ is satisfied, where $\Rb^\ast=\Ub^{-1/2}\Rb$. The column loading matrix $\Cb^\ast$ is normalized similarly, where $\Cb^\ast=\Vb^{-1/2}\Cb$. Let
\begin{equation}\nonumber
\vec(\Fb_t)\stackrel{i.i.d.}\sim\mathcal{N}(\textbf{0},\Db),\quad\textrm{where}\quad\Db=\textrm{diag}(1.5,1,0.5,1.5,1,0.5,1.5,1,0.5),
\end{equation}
so that the eigenvalues of $\bSigma_1$ in Assumption $\ref{asmpFactor}$ are distinct. Idiosyncratic matrices are generated according to $\vec(\Eb_t)\stackrel{i.i.d.}\sim\mathcal{MN}(\textbf{0},\Ub,\Vb)$. Thus, both $\{\Fb_t,1\leq t\leq T\}$ and $\{\Eb_t,1\leq t\leq T\}$ are independent across time, which simplifies the calculation of the asymptotic covariance matrix.
\begin{figure}[!ht]
  \centering
  \includegraphics[width=14cm, height=7cm]{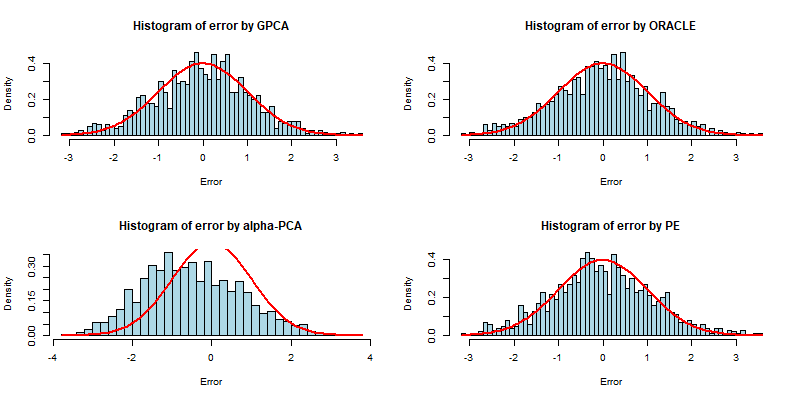}
 \caption{Histograms of estimation errors for $\Rb_{11}$ after normalization over 1000 replications under Case 1. $p_1=20$, $T=p_2=200$. Top left: GPCA. Top right: Oracle. Bottom left: $\alpha$-PCA with $\alpha=0$. Bottom right: PE. The red real line plots the probability density function of standard normal distribution.}\label{fig:1}
 \end{figure}
\begin{figure}[!ht]
  \centering
  \includegraphics[width=14cm, height=7cm]{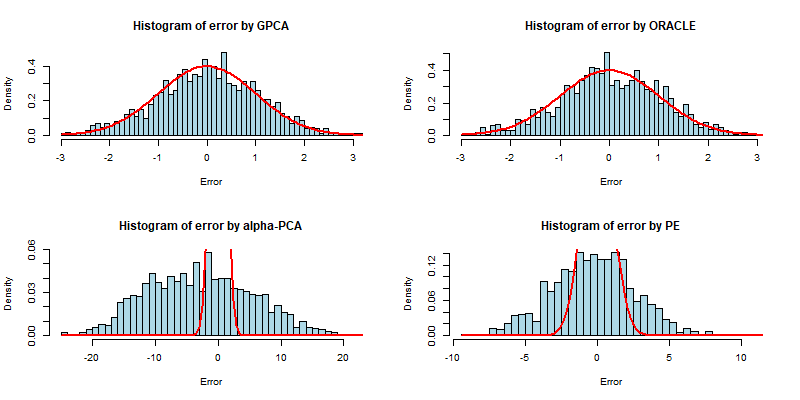}
 \caption{Histograms of estimation errors for $\Rb_{11}$ after normalization over 1000 replications under Case 2. $p_1=20$, $T=p_2=200$. Top left: GPCA. Top right: Oracle. Bottom left: $\alpha$-PCA with $\alpha=0$. Bottom right: PE. The red real line plots the probability density function of standard normal distribution.}\label{fig:2}
 \end{figure}

Figures $\ref{fig:1}$ and $\ref{fig:2}$ show the histograms of estimation errors for $\Rb_{11}$ after normalization over 1000 replications under Case 1 and 2, respectively. The asymptotic covariance matrices of the GPCA estimators and the Oracle estimators are the same theoretically, although the rotational matrices are not identical. First we focus on the case that $\Ub$ and $\Vb$ are weak cross-sectionally correlated
matrices. Figure $\ref{fig:1}$ shows that the GPCA, Oracle and PE estimators  are approximately normally distributed, while the estimator by $\alpha$-PCA deviates from normality. This result is expected as the condition $Tp_2=o(p_1^2)$ required by the  Theorem 2 in \cite{chen2023statistical}  is not yet met. Next, we consider the case where $\Ub$ and $\Vb$ are banded matrices. Figure $\ref{fig:2}$ shows that the GPCA, and Oracle estimators are approximately normally distributed, while the estimators by PE and $\alpha$-PCA exhibit significant deviations from normality. This further validates that our GPCA method can incorporate the additional information from $\Ub$ and $\Vb$, thereby improving estimation accuracy.

\section{Real Data Example}\label{Section5}
In this section, we demonstrate the empirical usefulness of our proposed data-driven GPCA method by analyzing a financial portfolio dataset ever studied by \cite{wang2019factor}, \cite{yu2022projected} and \cite{he2024matrix}. The financial portfolio dataset consists of monthly returns from 100 portfolios, which are well structured into a 10 × 10 matrix at each time point, with rows corresponding to 10 levels of market capital size (denoted as S1–S10) and columns corresponding to 10 levels of book-to-equity ratio (denoted as BE1-BE10). The dataset collects monthly returns from January 1964 to December 2019, covering a total of 672 months. Detailed information can be found on the website \url{http://mba.tuck.dartmouth.edu/pages/faculty/ken.french/data_library.html}.

We first apply the same preprocessing steps as in \cite{yu2022projected} and \cite{he2024matrix}, including standardization, missing value imputation, and stationarity test. The augmented Dickey-Fuller test rejects the null hypothesis, implying that all the series are stationary.
Additionally, we need to determine the number of factors $k_1$ and $k_2$. The IterER method proposed by \cite{yu2022projected} suggests that $(k_1,k_2)=(2,1)$,  the Rit-ER method proposed by \cite{he2024matrix} suggests that $(k_1,k_2)=(1,2)$ and the $\alpha$-PCA based ER method proposed by \cite{chen2023statistical} suggests $(k_1,k_2)=(1,1)$. For better illustration, we take $k_1=k_2=2$, also in view of  possible existence of weak factors. The estimated row and column loading matrices after varimax rotation and scaling are reported in Table $\ref{tab5}$.
\begin{table*}[!ht]
 	\begin{center}
 		\small
 		\addtolength{\tabcolsep}{1pt}
 		\caption{Loading matrices for Fama--French data set, after varimax rotation and scaling by 30. “GPCA”: data-driven GPCA method. “PE”: the projection estimation method by \cite{yu2022projected}. “$\alpha$-PCA”: the $\alpha$-PCA method in \cite{chen2023statistical} with $\alpha=0$. }\label{tab5}
 		 \renewcommand{\arraystretch}{1}
 		\scalebox{1}{ 		 \begin{tabular*}{14.5cm}{cc|cccccccccc}
 				\toprule[1.2pt]
 				 \multicolumn{12}{l}{Size}\\
 				\toprule[1.2pt]
 				 Method&Factor&S1&S2&S3&S4&S5&S6&S7&S8&S9&S10\\\hline
 \multirow{2}*{GPCA}&1&\cellcolor {Lavender}-31&\cellcolor {Lavender}-26&\cellcolor {Lavender}-22&\cellcolor {Lavender}-17&-13&-9&-4&-1&6&11
 	\\
 	&2&12&0&-10&-15&\cellcolor {Lavender}-23&\cellcolor {Lavender}-26&\cellcolor {Lavender}-29&\cellcolor {Lavender}-31&\cellcolor {Lavender}-34&\cellcolor {Lavender}-21
 	\\\hline
 				 \multirow{2}*{PE}&1&\cellcolor {Lavender}-16&\cellcolor {Lavender}-15&\cellcolor {Lavender}-12&\cellcolor{Lavender}-10&-8&-5&-3&-1&4&7\\
 				 &2&-6&-1&3&5&8&\cellcolor {Lavender}11&\cellcolor {Lavender}12&\cellcolor {Lavender}13&\cellcolor {Lavender}15&\cellcolor{Lavender}10\\\hline
 				 
 	 			 \multirow{2}*{$\alpha$-PCA}&1&\cellcolor {Lavender}-14&\cellcolor {Lavender}-14&\cellcolor {Lavender}-13&\cellcolor {Lavender}-11&-9&-7&-4&-2&3&7
 	\\
 	&2&-4&-2&1&3&6&9&\cellcolor {Lavender}12&\cellcolor {Lavender}13&\cellcolor {Lavender}16&\cellcolor {Lavender}14
 	\\

 				\bottomrule[1.2pt]		
\multicolumn{12}{l}{Book-to-Equity}\\
 				\toprule[1.2pt]
 				 Method&Factor&BE1&BE2&BE3&BE4&BE5&BE6&BE7&BE8&BE9&BE10\\\hline
 \multirow{2}*{GPCA}&1&11&1&-10&\cellcolor {Lavender}-16&\cellcolor {Lavender}-22&\cellcolor {Lavender}-26&\cellcolor {Lavender}-28&\cellcolor {Lavender}-29&\cellcolor {Lavender}-28&\cellcolor {Lavender}-24\\
 				 &2&\cellcolor {Lavender}43&\cellcolor {Lavender}34&\cellcolor {Lavender}24&\cellcolor {Lavender}16&9&4&0&-2&-1&1\\\hline
 				 \multirow{2}*{PE}&1&6&1&-4&-7&-10&\cellcolor {Lavender}-11&\cellcolor {Lavender}-12&\cellcolor {Lavender}-12&\cellcolor {Lavender}-12&-10
 				\\
 				&2&\cellcolor {Lavender}20&\cellcolor {Lavender}17&\cellcolor {Lavender}11&8&4&2&0&-1&-1&0\\\hline
 				 
 				 \multirow{2}*{$\alpha$-PCA}&1&6&2&-4&-7&-10&\cellcolor {Lavender}-11&\cellcolor {Lavender}-12&\cellcolor {Lavender}-13&\cellcolor {Lavender}-12&\cellcolor {Lavender}-11
\\
&2&\cellcolor {Lavender}19&\cellcolor {Lavender}18&\cellcolor {Lavender}12&8&4&2&0&-1&-1&-1
\\

 				\bottomrule[1.2pt]		
 		\end{tabular*}}		
 	\end{center}
 \end{table*}

From Table $\ref{tab5}$, we observe that the GPCA, PE and $\alpha$-PCA methods lead to very similar estimated loadings, where the GPCA method seems to lead to larger estimated loadings than the PE and $\alpha$-PCA methods. From the perspective of size, the small size portfolios load heavily on the first factor while the large size portfolios load mainly on the second factor. From the perspective of book-to-equity, the small BE portfolios load heavily on the second factor while the large BE portfolios load mainly on the first factor.

To further compare these methods, we use a rolling-validation procedure, as described in \cite{yu2022projected} and \cite{he2024matrix}. For each year $t$ from 1996 to 2019, we use the data from the previous $n$ years (the bandwidth) to fit the matrix factor model.  The fitted loadings are then used to estimate the factors and corresponding residuals for the 12 months of the current year. Specifically, let $\Yb_t^i$ and $\widetilde{\Yb}_t^i$ be the observed and estimated price matrix of month $i$ in year $t$, $\bar{\Yb}_t$ be the mean price matrix, and further define
\begin{equation}\nonumber
\textrm{MSE}_t=\frac{1}{12\times10\times10}\sum_{i=1}^{12}\left\|\widetilde{\Yb}_t^i-\Yb_t^i\right\|_F^2,\quad\rho_t=\frac{\sum_{i=1}^{12}\left\|\widetilde{\Yb}_t^i-\Yb_t^i\right\|_F^2}{\sum_{i=1}^{12}\left\|\Yb_t^i-\bar{\Yb}_t\right\|_F^2}, 
\end{equation}
as the mean squared pricing error and unexplained proportion of total variances, respectively. In the rolling-validation procedure, the variation of loading space is measured by $\upsilon_t:=\mathcal{D}(\widetilde{\Cb}_t\otimes\widetilde{\Rb}_t,\widetilde{\Cb}_{t-1}\otimes\widetilde{\Rb}_{t-1})$. 
\begin{table*}[!ht]
 	\begin{center}
 		\small
 		\addtolength{\tabcolsep}{1pt}
 		\caption{Rolling validation for the Fama-French portfolios. $12n$ is the sample size of the training set. $k_1=k_2=k$ is the number of factors. $\overline{\textrm{MSE}}$, $\bar{\rho}$, $\bar{\upsilon}$ are the mean pricing error, mean unexplained proportion of total variances and mean variation of the estimated loading space. “GPCA”: data-driven GPCA method. “PE”: the projection estimation method by \cite{yu2022projected}. “$\alpha$-PCA”: the $\alpha$-PCA method in \cite{chen2023statistical} with $\alpha=0$. }\label{tab6}
 		 \renewcommand{\arraystretch}{1}
 		\scalebox{1}{ 		 \begin{tabular*}{14.5cm}{ccccccccccc}
 				\toprule[1.2pt]
     \multirow{2}*{$n$}&\multirow{2}*{$k$}&\multicolumn{3}{l}{$\overline{\textrm{MSE}}$}&\multicolumn{3}{l}{$\bar{\rho}$}&\multicolumn{3}{l}{$\bar{\upsilon}$}\\
     \cmidrule(lr){3-5} \cmidrule(lr){6-8}\cmidrule(lr){9-11}
     
     &&GPCA&PE&$\alpha$-PCA&GPCA&PE&$\alpha$-PCA&GPCA&PE&$\alpha$-PCA
 				 \\
      \hline
 	5&1&\textbf{0.595}&0.870&0.862&\textbf{0.619}&0.802&0.796&\textbf{0.009}&0.176&0.241	\\
  10&1&\textbf{0.598}&0.855&0.860&\textbf{0.621}&0.784&0.791&\textbf{0.007}&0.085&0.203\\
  15&1&\textbf{0.598}&0.853&0.860&\textbf{0.621}&0.782&0.792&\textbf{0.009}&0.064&0.234\\
  \hline
  5&2&\textbf{0.595}&0.596&0.601&\textbf{0.620}&0.625&0.628&\textbf{0.009}&0.239&0.350\\
  10&2&\textbf{0.597}&0.601&0.611&\textbf{0.621}&0.628&0.636&\textbf{0.006}&0.092&0.261\\
  15&2&\textbf{0.598}&0.603&0.612&\textbf{0.621}&0.626&0.630&\textbf{0.008}&0.057&0.173\\
  \hline
 5&3&0.595&\textbf{0.522}&0.529&0.619&\textbf{0.550}&0.556&\textbf{0.008}&0.286&0.432\\
  10&3&0.598&\textbf{0.519}&0.526&0.621&\textbf{0.548}&0.555&\textbf{0.007}&0.114&0.353\\
  15&3&0.598&\textbf{0.517}&0.522&0.621&0.545&\textbf{0.544}&\textbf{0.008}&0.084&0.308\\
 				 			
 				\bottomrule[1.2pt]		
 		\end{tabular*}}		
 	\end{center}
 \end{table*}

Table $\ref{tab6}$ presents the mean values of MSE, $\rho$ and $\upsilon$ by the GPCA, PE and $\alpha$-PCA methods. Diverse combinations of bandwidth $n$ and the numbers of factors $(k_1=k_2=k)$ are compared. On the one hand, the pricing errors of the GPCA method are much lower than those of the PE and $\alpha$-PCA methods for $k=1$ and $k=2$. On the other hand, in terms of estimating the loading spaces, the GPCA method consistently performs more stably than the other methods. Financial data often exhibit complex heteroscedasticity, therefore the GPCA method is more preferred and acts as a robust approach for improving predictive accuracy.

\section{Discussion}\label{Section6}
In this study, we consider statistical inference for matrix factor model while focusing on heteroscedasticity of the idiosyncratic components. We start from the pseudo likelihood function and derive the pseudo likelihood estimators of the factor and loading spaces. {We derive the asymptotic distributions of the estimators of loadings, factors and common components when the covariance matrices are known in advance. In addition, we propose adaptive thresholding estimators for the separable covariance matrices and derive their convergence rates. We also show that this would not alter the asymptotic distributions of the Generalized Principal Component Analysis (GPCA)
estimators under some sparsity conditions.} 
Thorough numerical studies confirm the empirical advantage of the proposed GPCA method. We also demonstrate the practical usefulness of the proposed methods using a financial portfolio dataset. The limitation of our current work is that we do not provide theoretical guarantees  for the estimators along the solution path of the iterative algorithm, which encompasses both statistical and computational errors. We view it as an important direction for future research.

\section*{Acknowledgements}
He’s work is supported by National Science Foundation (NSF) of China (12171282), Taishan Scholars
Project Special Fund and Qilu Young Scholars Program of Shandong University. 

\bibliographystyle{model2-names}
\bibliography{Ref1}

\end{document}